\documentclass{amsart}
\usepackage{mathrsfs,url,amssymb,xypic}
\usepackage[francais,english]{babel}
\usepackage[latin1]{inputenc}

\newtheorem{thmm}{Theorem}[section]
\newtheorem{corr}[thmm]{Corollary}
\newtheorem{propp}[thmm]{Proposition}
\theoremstyle{definition}
\newtheorem{exee}[thmm]{Example}

\theoremstyle{plain}
\newtheorem{thm}{Theorem}[subsection]
\newtheorem{cor}[thm]{Corollary}
\newtheorem{lem}[thm]{Lemma}
\newtheorem{prop}[thm]{Proposition}
\theoremstyle{definition}
\newtheorem{defn}[thm]{Definition}
\newtheorem{cla}[thm]{Claim}
\newtheorem{rem}[thm]{Remark}
\newtheorem{que}[thm]{Question}
\newtheorem{conv}[thm]{Convention}
\newtheorem*{term}{Terminology}
\newtheorem{exe}[thm]{Example}
\numberwithin{equation}{section}

\newcommand{\eps}{\varepsilon}
\newcommand{\Z}{\mathbf{Z}}
\newcommand{\N}{\mathbf{N}}
\newcommand{\R}{\mathbf{R}}
\newcommand{\C}{\mathbf{C}}
\newcommand{\K}{\mathbf{K}}
\newcommand{\g}{\mathfrak{g}}
\newcommand{\mk}{\mathfrak}
\newcommand{\HH}{\mathscr{H}}
\newcommand{\SL}{\textnormal{SL}}
\newcommand{\SO}{\textnormal{SO}}
\newcommand{\SU}{\textnormal{SU}}
\begin{document}

%43A35 Positive definite functions on groups, semigroups, etc.
%20G25 Linear algebraic groups over local fields and their integers
%22D10 Unitary representations of locally compact groups
%22E40 Discrete subgroups of Lie groups
%22E50 Representations of Lie and linear algebraic groups over local fields

\selectlanguage{english}

\title[Relative Kazhdan Property]{Relative Kazhdan Property (Propriété~(T) relative)}
\author{Yves de Cornulier}%
\address{Yves de Cornulier\\
\'Ecole Polytechnique Fédérale de Lausanne (EPFL)\\
Institut de Géométrie, Algèbre et Topologie (IGAT)\\
CH-1015 Lausanne, Switzerland}

\email{decornul@clipper.ens.fr}

\date{\today}

% ----------------------------------------------------------------
\begin{abstract}
We perform a systematic investigation of Kazhdan's relative
Property~(T) for pairs $(G,X)$, where $G$ a locally compact group
and $X$ is any subset. When $G$ is a connected Lie group or a
$p$-adic algebraic group, we provide an explicit characterization
of subsets $X\subset G$ such that $(G,X)$ has relative
Property~(T). In order to extend this characterization to lattices
$\Gamma\subset G$, a notion of ``resolutions" is introduced, and
various characterizations of it are given. Special attention is
paid to subgroups of $\SU(2,1)$ and~$\SO(4,1)$.

\bigskip

\selectlanguage{french}

\noindent {\sc Résumé}.\hspace{1.5mm} Nous faisons une \'etude
syst\'ematique de la notion de propri\'et\'e (T) relative (de
Kazhdan) pour des paires $(G,X)$, o\`u $G$ est un groupe
localement compact et $X$ une partie quelconque. Lorsque $G$ est
un groupe de Lie connexe ou un groupe alg\'ebrique $p$-adique,
nous caract\'erisons de fa\c{c}on explicite les parties $X\subset
G$ telles que $(G,X)$ a la propri\'et\'e (T) relative. Une notion
convenable de ``r\'esolutions" permet d'\'etendre ces r\'esultats
aux r\'eseaux $\Gamma\subset G$, et nous en donnons diverses
caract\'erisations. Une attention particuli\`ere est port\'ee aux
sous-groupes de $\SU(2,1)$ et~$\SO(4,1)$.
\end{abstract}

\selectlanguage{english}

\maketitle
% ----------------------------------------------------------------

\section{Introduction}

Kazhdan's Property~(T) was introduced in a short paper by Kazhdan
\cite{K} in 1967. Since then, many consequences and
characterizations have been given by various authors.

The notion of {\em relative} Property for a pair $(G,N)$, where
$N$ is a normal subgroup in $G$ was implicit in Kazhdan's paper,
and later made explicit by Margulis \cite{Ma82}. The case when $H$
is an abelian normal subgroup is, by far, the best understood
\cite{K,Ma82,Bur,Sha99p,Sha99t}. However, it seems that it was
initially only considered as a technical tool. The most famous
case is the following: in order to prove Property~(T) for
$\SL_3(\R)$ (and other higher rank algebraic groups over local
fields), one uses, in most proofs, Property~(T) for the pair
$(\SL_2(\R)\ltimes\R^2,\R^2)$.

The definition of relative Property~(T) has been extended in
\cite{HV} to pairs $(G,H)$ with $H$ not necessarily normal in $G$.
Such pairs are extensively used in the work of Popa (see
\cite{Popa} and the references therein), in the context of
operator algebras. This motivated, for instance, new examples of
group pairs with relative Property~(T) of the form $(G\ltimes
N,N)$, with $N$ abelian \cite{Val04,Fer}.

We extend the definition of relative Property~(T) to pairs
$(G,X)$, where $X$ is any {\em subset} of $G$. The motivation for
this is that, given $G$, the knowledge of the family of subsets
$X$ such that $(G,X)$ has relative Property~(T) contains much
information about the unitary dual of $G$. It provides a knowledge
of $G$ much more precise than the bare information whether $G$ has
Property~(T). On the other hand, the family of subgroups with
relative Property~(T) provides a strictly weaker information (see
Example \ref{exe_intro_SO}).

Let $G$ be a locally compact group, and let $X\subset G$ be any
subset. We say that $(G,X)$ has relative Property~(T) if for every
net $(\varphi_i)$ of positive definite functions on $G$ that
converges to 1 uniformly on compact subsets, the convergence is
uniform on $X$.

In Section \ref{Sec_Prop_T} we establish various characterizations
of relative Property~(T) for a pair $(G,X)$, which were already
known \cite{Jol} in the case when $X$ is a subgroup. Here are the
main ones (the relevant definitions are recalled at the beginning
of Section \ref{Sec_Prop_T}).

\begin{thmm}[see Theorems \ref{thm:Del-Gui_relT},
\ref{thm:propT_pure_phi}, and \ref{thm:T_rep_irr}]
Let $G$ be a locally compact, $\sigma$-compact group, and
$X\subset G$ any subset. The following are equivalent:
\begin{enumerate}
\item\label{item:i1} {\em(Positive definite functions)} $(G,X)$
has relative Property~(T).

\item {\em(Representations with almost invariant vectors)} For
every $\eps>0$, and for every unitary representation $\pi$ of $G$
that almost has invariant vectors, $\pi$ has $(X,\eps)$-invariant
vectors.

\item {\em(Kazhdan pairs)} For every $\eps>0$, there exists a
compact subset $K\subset G$ and $\eta>0$ such that, every unitary
representation of $G$ that has a $(K,\eta)$-invariant vector has a
$(X,\eps)$-invariant vector.

\item {\em(Conditionally negative definite functions)} Every
conditionally negative definite function on $G$ is bounded on $X$.

\item\label{item:prop_FH_rel} {\em(Isometric actions on affine
Hilbert spaces)} For every affine, isometric action of $G$ on a
affine Hilbert space $\mathscr{H}$, and every $v\in\mathscr{H}$,
$Xv$ is bounded.

\item {\em(Topology of the unitary dual)} For every $\eps>0$, and
for every net $(\pi_i)$ of irreducible unitary representations of
$G$ that converges to $1_G$, eventually $\pi_i$ has a
$(X,\eps)$-invariant vector.
\end{enumerate}\label{thm:intro_relT_eq}
\end{thmm}

Recall that a locally compact group $G$ is Haagerup if it has a
net of $C_0$ positive definite functions that converges to 1,
uniformly on compact subsets. It is clear from the definition that
if $G$ is Haagerup, then, for every $X\subset G$, the pair $(G,X)$
has relative Property~(T) if and only if $X$ is relatively compact
in $G$. The question whether the converse holds was asked (in a
slightly different formulation) in \cite{AW2}. We say that $G$
satisfies the TH alternative if it is either Haagerup, or has a
subset $X$ with noncompact closure, such that $(G,X)$ has relative
Property~(T). The question becomes: does there exist a locally
compact group that does not satisfy the TH alternative? We leave
it open.

\medskip

By a result of Kazhdan, if $G$ is a locally compact group with
Property~(T), then $G$ is compactly generated. The same argument
shows that if $(G,X)$ has Property~(T), then $X$ is contained in a
compactly generated subgroup of $G$. Here is a stronger result,
which says, in a certain sense, that all questions about relative
Property~(T) reduce to the compactly generated case.

\begin{thmm}[Theorem \ref{thm:propT_and_cpt_gen}]
Let $G$ be a locally compact group, and $X\subset G$. Then $(G,X)$
has relative Property~(T) if and only if there exists an open,
compactly generated subgroup $H$ of $G$, containing $X$, such that
$(H,X)$ has relative Property~(T).\label{thm:intro_cpt_gen}
\end{thmm}

Note that Theorem \ref{thm:intro_cpt_gen} is new even in the case
when $X\subset G$ is a normal subgroup. As a corollary of Theorem
\ref{thm:intro_cpt_gen} (see Remark \ref{r:dir_lim_RL_Haag}), we
deduce that a locally compact group satisfies the TH alternative
if and only if all its open, compactly generated subgroup do.

We are interested in the question of determining, given a group
$G$, subsets $X$ such that $(G,X)$ has relative Property~(T). As a
general result, we show, provided that $G$ is compactly generated,
that such subsets coincide with the bounded subsets for a
well-defined, essentially left-invariant metric on $G$, which we
call the $H$-metric (see subsection \ref{subs:Hmetric}).

\medskip

In Section \ref{Sec_relT_lie_alg}, we focus on relative
Property~(T) in connected Lie groups and linear algebraic groups
over a local field $\K$ of characteristic zero.

Let $G$ be a connected Lie group. Let $R$ be its radical, and $S$
a Levi factor. Define $S_{nc}$ as the sum of all non-compact
factors of $S$, and $S_{nh}$ as the sum of all factors of $S_{nc}$
with Property~(T). Finally define the T-radical
$R_T=\overline{S_{nh}[S_{nc},R]}$. It is easily checked to be a
characteristic subgroup of $G$.

\begin{thmm}[Theorem \ref{thm:Trel_Lie}]
$(G,R_T)$ has relative Property~(T).\label{thm:intro_relT_lie}
\end{thmm}

On the other hand, by results in \cite[Chap. 4]{CCJJV}, $G/R_T$
has the Haagerup Property. As a consequence:

\begin{corr}
Let $X\subset G$, and $p$ denote the projection $G\to G/R_T$. Then
$(G,X)$ has relative Property~(T) if and only if $\overline{p(X)}$
is compact.\label{cor:intro:relTlie}
\end{corr}

The proof of Theorem \ref{thm:intro_relT_lie} makes use of the
following proposition:

\begin{propp}
Let $S$ be a semisimple connected Lie group without compact
factors, and $R$ a nilpotent, simply connected Lie group, endowed
with an action of $S$, and set $G=S\ltimes R$. Suppose that
$[S,R]=R$. Then $(G,R)$ has relative Property~(T).
\end{propp}

The proposition is proved as follows: we work by induction on
$\dim(R)$, we pick a nontrivial central subgroup $V$ of $R$,
normal in $G$, and we can reduce to two cases. Either $[S,V]=V$,
so that, by well-known results (which can be attributed to
Kazhdan, Margulis, Burger), the pair $(G,V)$ has relative
Property~(T), and the result follows by induction, or $V$ is
central in $G$. To handle this case, we formulate an ad-hoc result
of stability of relative Property~(T) by central extensions.

The case of a linear algebraic group $G$ over a local field $\K$
of characteristic zero is similar. Let $R$ be the radical, $S$ a
Levi factor, and define $S_{nc}$ and $S_{nh}$ as in the case of
Lie groups, and set $R_T=S_{nh}[S_{nc},R]$; this is a closed
characteristic subgroup.

\begin{thmm}[Theorem \ref{thm:RT_alg}]
$(G(\K),R_T(\K))$ has relative
Property~(T).\label{thm:intro_relT_alg}
\end{thmm}

On the other hand, it is easily checked that $G/R_T$ has the
Haagerup Property. A corollary similar to Corollary
\ref{cor:intro:relTlie} follows.

Another corollary is the following result, already known:

\begin{corr}
Let $G$ be a connected Lie group (respectively a linear algebraic
group over~$\K$).

\begin{itemize}
    \item[(1)] $G$ [respectively $G(\K)$] has Property~(T) if and only
    if $R_T$ [resp. $R_T(\K)$] is cocompact in $G$ [resp. in $G(\K)$].

    \item[(2)] $G$ [resp. $G(\K)$] is Haagerup if and only if $R_T=1$.
\end{itemize}
\end{corr}
Assertion (1) is a result of S.P. Wang \cite{Wang}, and (2) is due
to \cite[Chap. 4]{CCJJV} for connected Lie groups and to
\cite{CorJLT} in the $p$-adic case.

\medskip

Section \ref{Sec_resol} is devoted to explain how these phenomena
are inherited by subgroups of finite covolume.

Let $G$ be a locally compact group, $N$ a closed, normal subgroup,
and $H$ a subgroup of finite covolume in $G$. It is known (see
\cite{Jol}) that if $(G,N)$ has relative Property~(T), then so
does $(H,N\cap H)$. However, this result is of limited use insofar
as $N\cap H$ may be small (for instance, reduced to $\{1\}$) even
if $N$ is noncompact: this phenomenon is very frequent in the
context of irreducible lattices in products of algebraic groups
over local fields.

We need a definition that enlarges the notion of relative
Property~(T) of normal subgroups. The datum of $N$ normal in $G$
is equivalent to the datum of the morphism $G\to G/N$. More
generally, we are led to consider arbitrary locally compact groups
$G,Q$, and a morphism $f:G\to Q$ with dense image. We say that $f$
is a {\it resolution} if, for every unitary representation $\pi$
of $G$ almost having invariant vectors, $\pi$ has a nonzero
subrepresentation $\rho$ factoring through a representation
$\tilde{\rho}$ of $Q$, and $\tilde{\rho}$ almost has invariant
vectors.

Given a closed normal subgroup $N$ in a locally compact group $G$,
$(G,N)$ has relative Property~(T) if and only if $G\to G/N$ is a
resolution. In view of Theorems \ref{thm:intro_relT_lie} and
\ref{thm:intro_relT_alg}, a wealth of examples of resolutions are
provided by the following result, essentially due to Margulis
\cite[Chap. III, Section 6]{Margulis}, and which also uses
arguments borrowed from \cite{BL}.

\begin{thmm}[Theorem \ref{thm:MBL}]
Let $G$ be a locally compact group, $N$ a closed, normal subgroup.
Suppose that $(G,N)$ has relative Property~(T) (equivalently, the
projection $p:G\to G/N$ is a resolution).

Let $H$ be a closed subgroup of finite covolume in $G$, and write
$Q=\overline{p(H)}$. Then $p:H\to Q$ is a
resolution.\label{thm:intro_MBL}
\end{thmm}

Resolutions allow to prove compact generation of some locally
compact groups. The following theorem generalizes Proposition 2.8
of \cite{LZ}.

\begin{thmm}[Theorem \ref{thm:resol_cpt_gen}]
Let $G\to Q$ be a resolution. Then $G$ is compactly generated if
and only if $Q$ is.
\end{thmm}

Thus, compact generation can be said to be ``invariant under
resolutions". We provide some other examples.

\begin{propp}[see Theorem \ref{thm:resol_cpt_gen}, Corollary
\ref{cor:resol_FA}, Proposition \ref{prop:resol_FAb_etc}]Let $G\to
Q$ be a resolution. Then, if (P) is one of the properties below,
then $G$ has Property (P) if and only if $Q$ does:

\begin{itemize}
    \item Property~(T),
    \item Compact generation,
    \item Property (FA):
every isometric action on a tree has a fixed point,
    \item Every isometric
action on a Euclidean space has a fixed point,
    \item Property
($\tau$): the trivial representation is isolated among irreducible
unitary representations with finite image.
\end{itemize}
\end{propp}

Resolutions give rise to pairs with relative Property~(T).

\begin{propp}
Let $p:H\to Q$ be a resolution. Given any subset $X\subset H$, if
$(Q,p(X))$ has relative Property~(T), then so does $(H,X)$. In
particular, if $Q$ is Haagerup, then $(H,X)$ has relative
Property~(T) if and only if $\overline{p(X)}$ is
compact.\label{prop:intro_resol_relT}
\end{propp}

The first consequence is that if $Q$ is Haagerup, then $H$
satisfies the TH alternative (Corollary
\ref{cor:resol_Haag_TH_alt}). In view of Theorems
\ref{thm:intro_relT_alg} and \ref{thm:intro_MBL}, this applies to
lattices in products of real and $p$-adic algebraic groups.

We can also derive some new phenomena of relative Property~(T).

\begin{exee}[see Proposition \ref{prop:SO5_ltimes_R5},
the proof of Proposition \ref{prop:sg_trel_SO41_SU21}(2), and
Remark \ref{rem:SO3_ltimes_R3}]

Consider the group
$\Gamma_n=\SO_n(\Z[2^{1/3}])\ltimes\Z[2^{1/3}]^n$.

If $n\le 2$, it is solvable, hence Haagerup. If $n=3,4$, it was
observed in \cite{CorJLT}, that $\Gamma_n$ is not Haagerup, but
has no infinite subgroup with relative Property~(T). Thanks to
resolutions, we can see this more concretely: the natural morphism
$i:\Gamma_n\to \SO_n(\Z[2^{1/3}])\ltimes\R^n$ is a resolution. It
follows that $(\Gamma_n,B_n)$ has relative Property~(T), where
$B_n$ is the intersection of the unit ball of $\R^n$ with
$\Z[2^{1/3}]^n$.

If $n\ge 5$, $\Gamma_n$ has a different behaviour, due to the fact
that $\SO_n(\Z[2^{1/3}])$ has Property~(T) (whereas it is Haagerup
if $n\le 4$). It follows that, if $n\ge 5$, then
$\Gamma_n\to\SO_n(\R)\ltimes\R^n$ is a resolution. We deduce an
interesting property for $\Gamma_n$: it has Property $(\tau)$, but
has a finite dimensional unitary representation $\pi$ such that
$H^1(\Gamma_n,\pi)\neq 0$.\label{exe_intro_SO}
\end{exee}

As pointed out in \cite{PP}, all previously known examples of
group pairs with relative Property~(T) were derived from group
pairs where the subgroup is normal. For instance, for every group
$G$, $((\SL_2(\Z)\ltimes\Z^2)\ast G,\Z^2)$ has relative
Property~(T); the only role of $G$ is to prevent $\Z^2$ from being
normal.

Resolutions and Proposition \ref{prop:intro_resol_relT} allows us
to find group pairs with relative Property~(T) that do not derive
in such a way from group pairs with a normal subgroup. This is
illustrated by projections of irreducible lattices from
$\SO(4,1)\times\SO(5,\C)$, such as $\SO(4,1)(\Z[2^{1/3}])$, see
Theorem \ref{thm:sgSO41SU21} and the proof of Proposition
\ref{prop:sg_trel_SO41_SU21}(1).

\begin{propp}
Let $G$ be either $\SO(4,1)$ or $\SU(2,1)$, and $\Gamma$ a
subgroup, viewed as a discrete group (but not necessarily discrete
in $G$).

\begin{itemize}
    \item[1)] If $\Lambda\subset\Gamma$ is a normal subgroup and
$(\Gamma,\Lambda)$ has relative Property~(T), then $\Lambda$ is
finite.
     \item[2)] There exists $\Gamma$ containing an infinite
subgroup $\Lambda$ such that $(\Gamma,\Lambda)$ has relative
Property~(T).\end{itemize}
\end{propp}

Finally, we prove various equivalences for resolutions. They
generalize known equivalences for pairs $(G,N)$ with $N$ normal in
$G$, but the proofs are more involved.

\begin{thmm}[see Theorems \ref{thm:resol} and \ref{thm:resol_irr}]
Let $G,Q$ be a locally compact $\sigma$-compact groups, and let
$f:G\to Q$ be a morphism with dense image. The following are
equivalent:

\begin{enumerate}

\item $f:G\to Q$ is a resolution.

\item $f$ is an ``affine resolution": for every isometric affine
action of $G$ on a Hilbert space, there exists a $G$-invariant
nonempty closed subspace on which the action factors through $Q$.

\item For every subset $X\subset G$ such that $\overline{f(X)}$ is
compact, $(G,X)$ has relative Property~(T).

\item For every net $(\pi_i)$ of irreducible unitary
representations of $G$, if $\pi_i\to 1_G$, then eventually $\pi_i$
factors through a representation $\tilde{\pi}_i$ of $Q$, and
$\tilde{\pi}_i\to 1_Q$.
\end{enumerate}
\end{thmm}

\bigskip

\noindent\textbf{Acknowledgments.} I thank Bachir Bekka, Pierre de
la Harpe, Vincent Lafforgue, Romain Tessera, and Alain Valette for
valuable comments and useful discussions.

\section{Property~(T) relative to subsets}\label{Sec_Prop_T}

\begin{term}
Throughout the paper, by {\it morphism} between topological groups
we mean continuous group homomorphisms. If $X$ is a Hausdorff
topological space, a subset $Y\subset X$ is {\it relatively
compact} if its closure in $X$ is compact.

If $G$ is a group, a {\em positive definite} function on $G$ is a
function: $\varphi:G\to\C$ such that
$\varphi(g)=\overline{\varphi(g^{-1})}$ for all $g\in G$, and, for
all $n\in\N$, $(g_1,\dots,g_n)\in G^n$ and
$(c_1,\dots,c_n)\in\C^n$, the inequality
$\sum_{i,j}c_i\overline{c_j}\varphi(g_i^{-1}g_j)\ge 0$ holds. We
say that a positive definite function $\varphi$ is {\em
normalized} if $\varphi(1)=1$.

If $G$ is a group, a {\em conditionally negative definite}
function on $G$ is a function: $\psi:G\to\R$ such that
$\psi(1)=0$, $\psi(g)=\psi(g^{-1})$ for all $g\in G$, and for all
$n\in\N$, $(g_1,\dots,g_n)\in G^n$ and $(c_1,\dots,c_n)\in\R^n$
such that $\sum c_i=0$, the inequality
$\sum_{i,j}c_ic_j\psi(g_i^{-1}g_j)\le 0$ holds.

A locally compact group $G$ is \textit{Haagerup} if for every
compact $K\subset G$ and $\eps>0$, there exists a real-valued
normalized positive definite function in $G$, vanishing at
infinity, and such that $\varphi|_K>1-\eps$. A locally compact
group $G$ is \textit{a-T-menable} if it has a proper conditionally
negative definite function. These two properties are introduced
and called (3A) and (3B) in \cite{AW2}, where the authors prove
their equivalence in the $\sigma$-compact case.

All unitary representations $\pi$ of a locally compact group are
supposed continuous, that is, the function $g\mapsto\pi(g)\xi$ is
continuous for every $\xi$ in the representation space. Recall
that $1_G\prec\pi$ (or $1\prec G$ if there is no ambiguity) means
that the unitary representation $\pi$ \textit{almost has invariant
vectors}, that is, for every compact subset $K\subset G$ and every
$\eps>0$, there exists a nonzero vector $\xi$ such that
$\sup_{g\in K}\|\pi(g)\xi-\xi\|\le\eps\|\xi\|$.
\end{term}

\subsection{Property~(T) relative to subsets}

\begin{defn}
Let $G$ be a locally compact group, and $X$ any subset. We say
that $(G,X)$ has relative Property~(T) if, for every net
$(\varphi_i)$ of continuous normalized positive definite functions
that converges to 1 uniformly on compact subsets, the convergence
is uniform on $X$.

We say that $(G,X)$ has relative Property (FH) if every continuous
conditionally definite negative function on $G$ is bounded on $X$.
\end{defn}

If $\varphi$ is a positive definite function on $G$, then so is
$|\varphi|^2$. Thus, the definition of relative Property~(T)
remains unchanged if we only consider real-valued positive
definite functions or even non-negative real-valued positive
definite functions.

\begin{que}[\cite{AW2}]
For a locally compact group $G$, consider the two following
properties.
\begin{itemize}\item[(1)] $G$ is a-T-menable;

\item[(2)] for every subset $X\subset G$, the pair $(G,X)$ has
relative Property (FH) if and only if $\overline{X}$ is
compact.\end{itemize}

(Note that the implication (1)$\Rightarrow$(2) is trivial.) Does
there exist a $\sigma$-compact, locally compact group $G$
satisfying (2) and not (1)? \label{que:Haag=nonrelT}
\end{que}

\begin{rem}
If $G$ is locally compact but not $\sigma$-compact, (1) and (2) of
Question \ref{que:Haag=nonrelT} are not equivalent: (1) is always
false, while a characterization of (2) is less clear. For
instance, if $G$ is any abelian locally compact group, then (2) is
fulfilled. On the other hand, it was shown in \cite{CorUFH} that
if $F$ is a non-nilpotent finite group, then $F^\N$, viewed as a
discrete group, does not satisfy (2); moreover, if $F$ is perfect,
then $F^\N$ has Property (FH). Besides, being locally finite,
these groups are amenable, hence Haagerup.
\end{rem}

It is maybe worth comparing Question \ref{que:Haag=nonrelT} to the
following result, essentially due to \cite{GHW}:

\begin{prop}
The locally compact, $\sigma$-compact group $G$ is a-T-menable if
and only if there exists a sequence $(\psi_n)_{n\in\N}$ of
continuous conditionally negative definite functions on $G$, such
that, for every sequence $(M_n)$ of positive numbers, $\{g\in G|\;
\forall n,\; \psi_n(g)\le M_n\}$ is compact.
\end{prop}
\begin{proof} The direct implication is trivial (take any proper
function $\psi$, and $\psi_n=\psi$ for all $n$). Conversely,
suppose the existence of a family $(\psi_n)$ satisfying the
condition. Let $(K_n)$ be an increasing sequence of compact
subsets of $G$ whose interiors cover $G$. There exists a sequence
$(\eps_n)$ such that $\eps_n\psi_n\le 2^{-n}$ on $K_n$. Set
$\psi=\sum_n\eps_n\psi_n$; since the series is convergent
uniformly on compact subsets, $\psi$ is well-defined and
continuous. Then, for every $M<\infty$, the set $\{\psi\le M\}$ is
contained in $\{g|\;\forall n,\; \psi_n\le M/\eps_n\}$, which is,
by assumption, compact.\end{proof}

\begin{que}%[\cite{AW2}]
Does there exist a locally compact group $G$ without the Haagerup
Property, but with no unbounded subset $X$ such that $(G,X)$ has
relative Property~(T)?\label{que:Haag=wdRL?}
\end{que}

\begin{rem}
Akemann and Walter \cite{AW1} introduce the following definition:
a locally compact group has the weak dual Riemann-Lebesgue
Property if, for every $\eps,\eta>0$ and every compact subset $K$
of $G$, there exists a compact subset $\Omega$ of $G$ such that,
for every $x\in G-\Omega$, there exists a normalized, real-valued,
positive definite function $\varphi$ on $G$ such that
$\varphi(x)\le\eta$ and $\sup\{1-\varphi(g)|\;g\in K\}\le\eps$.

%(where $\|f\|_{\infty}^K$ denotes $\sup_{x\in K}|f(x)|$).
%
It can be shown \cite[Proposition 2.7]{CorThe} that a locally
compact group $G$ has the weak-dual Riemann-Lebesgue if and only
if every subset $X\subset G$ such that $(G,X)$ has relative
Property~(T) is relatively compact. Akemann and Walker \cite{AW2}
ask if the weak dual Riemann-Lebesgue Property is equivalent to
the Haagerup Property; this question is therefore equivalent to
Question \ref{que:Haag=wdRL?}.
\end{rem}

\begin{rem}
1) It follows from Theorem \ref{thm:Del-Gui_relT} that Questions
\ref{que:Haag=nonrelT} and \ref{que:Haag=wdRL?} are equivalent for
locally compact, $\sigma$-compact groups.

2) It follows from \cite[Proposition 6.1.1]{CCJJV} and Theorem
\ref{thm:propT_and_cpt_gen} that if Question \ref{que:Haag=wdRL?}
has a positive answer, then the example can be chosen compactly
generated.\label{r:dir_lim_RL_Haag}
\end{rem}

\begin{defn}We say that $G$ satisfies the TH alternative if it is
either Haagerup, or has a subset $X$ with noncompact closure, such
that $(G,X)$ has relative Property~(T).
\end{defn}

Question \ref{que:Haag=wdRL?} becomes: does there exist a locally
compact group not satisfying the TH alternative?

\begin{rem}
Here is an obstruction to the Haagerup Property for a locally
compact, compactly generated group $G$, that does not formally
imply the existence of a non-relatively compact subset with
relative Property~(T). Let $\omega$ belong to the Stone-\v Cech
boundary $\beta G\smallsetminus G$ of $G$. Let us say that
$(G,\omega)$ has relative Property~(T) if, for every conditionally
negative definite function $\psi$ on $G$, its canonical extension
$\tilde{\psi}:\beta G\to\R_+\cup\{\infty\}$ satisfies
$\tilde{\psi}(\omega)<\infty$.

It is clear that relative Property~(T) for $(G,\omega)$ prevents
$G$ from being Haagerup. On the other hand, I see no reason why
this should imply the existence of a non-relatively compact subset
with relative Property~(T).
\end{rem}

\subsection{Various equivalences}

\begin{defn}
Let $G$ be a locally compact group and $X\subset G$. Given a
unitary representation $\pi$ of $G$ and $\eps\ge 0$, a
$(X,\eps)$-invariant vector for $\pi$ is a nonzero vector $\xi$ in
the representation space such that
$\|\pi(g)\xi-\xi\|\le\eps\|\xi\|$ for every $g\in X$.
\end{defn}

\begin{defn}
Let $G$ be a locally compact group, $X,W$ subsets, and
$\eps,\eta>0$. We say that $(W,\eta)$ is a $\eps$-Kazhdan pair for
$(G,X)$ if, for every unitary representation $\pi$ of $G$ that has
a $(W,\eta)$-invariant vector, $\pi$ has a $(X,\eps)$-invariant
vector. Given $G,X,W,\eps$, if such $\eta>0$ exists, we say that
$W$ is a $\eps$-Kazhdan subset for $(G,X)$.
\end{defn}

The following result generalizes a result due to Jolissaint
\cite{Jol} when $X$ is a subgroup.

\begin{thm}
Let $G$ be a locally compact group, and let $X\subset G$ be a
subset. Consider the following properties.
\begin{itemize}

\item[(1)] $(G,X)$ has relative Property~(T).

\item[(2)] For every $\eps>0$, there exists a compact
$\eps$-Kazhdan subset for $(G,X)$.

\item[(2')] For some $\eps<\sqrt{2}$, there exists a compact
$\eps$-Kazhdan subset for $(G,X)$.

\item[(3)] For every $\eps>0$ and every unitary representation
$\pi$ of $G$ such that $1\prec\pi$, the representation $\pi$ has a
$(X,\eps)$-invariant vector.

\item[(3')] There exists $\eps<\sqrt{2}$ such that for every
unitary representation $\pi$ of $G$ satisfying $1\prec\pi$, the
representation $\pi$ has a $(X,\eps)$-invariant vector.

\item[(4)] $(G,X)$ has relative Property (FH), i.e. satisfies
(\ref{item:prop_FH_rel}) of Theorem \ref{thm:intro_relT_eq}.

\end{itemize}

Then the following implications hold:
$$\xymatrix{ (1) \ar@{=>}[r] & (2)  \ar@{=>}[d] \ar@{=>}[r] & (3)  \ar@{=>}[d] & \\
                             & (2') \ar@{=>}[r]             & (3') \ar@{=>}[r] & (4). }$$

Moreover, if $G$ is $\sigma$-compact, then (4)$\Rightarrow$(1), so
that they are all equivalent. \label{thm:Del-Gui_relT}
\end{thm}

\begin{proof} (1)$\Rightarrow$(2) Suppose the contrary. There exists
$\eps>0$ such that, for every $\eta>0$ and every compact subset
$K\subset G$, there exists a unitary representation $\pi_{\eta,K}$
of $G$ that has a $(K,\eta)$-invariant unit vector $\xi_{\eta,K}$,
but has no $(X,\eps)$-invariant vector. Denote by
$\varphi_{\eta,K}$ the corresponding coefficient. Then, when
$\eta\to 0$ and $K$ becomes big, $\varphi_{\eta,K}$ converges to
1, uniformly on compact subsets. By relative Property~(T), the
convergence is uniform on $X$. It follows that, for some $K$ and
some $\eta$, the representation $\pi_{\eta,K}$ has a
$\eps$-invariant vector, a contradiction.

(2)$\Rightarrow$(2'), (2)$\Rightarrow$(3), (2')$\Rightarrow$(3'),
and (3)$\Rightarrow$(3') are immediate.

(3')$\Rightarrow$(4) Let $\psi$ be a conditionally negative
definite function on $G$, and, for $t>0$, let
$(\pi_t,\mathcal{H}_t)$ be the cyclic unitary representation of
$G$ associated with the positive definite function $e^{-t\psi}$.
Set $\rho_t=\pi_t\otimes\overline{\pi_t}$. Since $\pi_t\to 1_G$
when $t\to 0$, so does $\rho_t$.

Suppose that $\psi$ is not bounded on $X$: $\psi(x_n)\to\infty$
for some sequence $(x_n)$ in $X$. Then we claim that for every
$t>0$ and every $\xi\in\mathcal{H}_t\otimes
\overline{\mathcal{H}_t}$, we have
$\langle\rho_t(x_n)\xi,\xi\rangle\to 0$ when $n\to\infty$.
Equivalently, for every $\xi\in\mathcal{H}_t\otimes
\overline{\mathcal{H}_t}$ of norm one, $\|\rho_t(x_n)\xi-\xi\|\to
\sqrt{2}$. This is actually established in the proof of
\cite[Lemma 2.1]{Jol} (where the assumption that $X=H$ is a
subgroup is not used for this statement).

By Lebesgue's dominated convergence Theorem, it follows that if
$\rho$ denotes the representation $\bigoplus_{t>0}\rho_t$, then
$\langle\rho(x_n)\xi,\xi\rangle\to 0$ for every $\xi$. In
particular, for every $\eps<\sqrt{2}$, the representation $\rho$
has no $(X,\eps)$-invariant vector. Since $1\prec\rho$, this
contradicts~(3').

(4)$\Rightarrow$(1) The proof is a direct adaptation of that of
the analogous implication in \cite[Theorem~3]{AW2} and we do not
repeat it here.\end{proof}

\begin{rem}
When $X=H$ is a subgroup, we retrieve a result of \cite{Jol}. Note
that, in this case, by a well-known application of the ``Lemma of
the centre" \cite[Lemma 2.2.7]{BHV}, Condition (2') of Theorem
\ref{thm:Del-Gui_relT} can be chosen with $\eps=0$, i.e. becomes:
for every unitary representation of $G$ such that $1\prec\pi$,
there exists a nonzero $H$-invariant vector.
\end{rem}

\begin{rem}
When $G$ is not $\sigma$-compact, whether the implication
(3')$\Rightarrow$(1) holds is not known, except when $X$ is a
normal subgroup \cite{Jol}. On the other hand, it is known
\cite{CorUFH} that (4)$\Rightarrow$(3') does not always hold for
general locally compact groups, even when we assume that $X=G$ is
a discrete group.
\end{rem}

\subsection{Relative Property~(T) can be read on irreducible unitary representations.}

The following lemma, due to Choquet (unpublished), is proved in
\cite[B.14 p. 355]{Dix}.

\begin{lem}
Let $K$ be a compact, convex subset of a locally convex space $E$.
Let $x$ be an extremal point of $K$. Let $\mathcal{W}$ be the set
of all open half-spaces of $E$ that contain $x$. Then $\{W\cap
K|\; W\in\mathcal{W}\}$ is a neighbourhood basis of $x$ in
$K$.\qed\label{lem:cvx_extremal_basis}
\end{lem}

Denote $\mathcal{P}(G)$ [resp. $\mathcal{P}_1(G)$, resp.
$\mathcal{P}_{\le 1}(G)$] the set of all (complex-valued) positive
definite function $\varphi$ on $G$ [resp. such that
$\varphi(1)=1$, resp. such that $\varphi(1)\le 1$].

Recall that $\varphi\in\mathcal{P}(G)$ is {\it pure} if it
satisfies one of the two equivalent conditions: (i) $\varphi$ is
associated to an irreducible unitary representation; (ii)
$\varphi$ belongs to an extremal axis of the convex cone
$\mathcal{P}(G)$.

\begin{thm}
Let $G$ be a locally compact group and $X$ a subset. The following
are equivalent:
\begin{itemize}\item[(i)] $(G,X)$ has relative Property~(T).

\item[(ii)] For every net of continuous, normalized {\em pure}
positive definite functions on $G$ converging to 1, the
convergence is uniform on
$X$.\end{itemize}\label{thm:propT_pure_phi}
\end{thm}
\begin{proof} (i)$\Rightarrow$(ii) is trivial; conversely suppose that $G$ satisfies
(ii). Endow the space $L^\infty(G)=L^1(G)^*$ with the weak*
topology. Let $\mathcal{W}$ be the set of all open half-spaces of
$L^\infty(G)$ containing the constant function $1$. Finally set
$K=\mathcal{P}_{\le 1}(G)$.

Recall Raikov's Theorem \cite[Th\'eor\`eme 13.5.2]{Dix}: on
$\mathcal{P}_1(G)$, the weak* topology coincides with the topology
of uniform convergence on compact subsets. Let $\mathcal{L}$ be
the set of all continuous linear forms $u$ on $L^\infty(G)$ such
that $u(1)=1$ and $u|_K\le 1$. Since $K$ is convex and compact for
the weak*-topology, by Lemma \ref{lem:cvx_extremal_basis},
$\{\{u>1-\eps\}\cap K|\; u\in\mathcal{L},\eps>0\}$ is a basis of
open neighbourhoods of $1$ in $\mathcal{P}_{\le 1}(G)$.

Hence, by (ii), and using Raikov's Theorem, for every $1>\eps>0$,
there exists $u\in\mathcal{L}$ and $\eta>0$ such that, for every
pure $\varphi\in\mathcal{P}_1(G)$, $u(\varphi)>1-\eta$ implies
$\varphi\ge 1-\eps$ on $X$.

Let $\varphi=\sum \lambda_i\varphi_i$ be a convex combination of
continuous, normalized, pure positive definite functions
$\varphi_i$. Suppose that $u(\varphi)>1-\eta\eps$. Decompose
$\varphi$ as $\sum \lambda_j\varphi_j+\sum \lambda_k\varphi_k$,
where $u(\varphi_j)>1-\eta$ and $u(\varphi_k)\le 1-\eta$. Then
$$1-\eta\eps\;<
\;u(\varphi)\;= \;\sum \lambda_ju(\varphi_j)+\sum
\lambda_ku(\varphi_k)$$ $$\le\; \sum \lambda_j+\sum
\lambda_k(1-\eta)\;=\;1-\eta\sum\lambda_k,$$ so that
$\sum\lambda_k\le\eps$. Hence, on $X$, we have
$$\varphi\;=\;\sum\lambda_j\varphi_j+\sum\lambda_k\varphi_k\;\ge\;
\sum\lambda_j(1-\eps)-\sum\lambda_k$$ $$\ge\;
\sum\lambda_j-2\eps\;\ge\; 1-3\eps.$$

Set $K_{u,\eps\eta}=\{\varphi\in K|\;u(\varphi)>1-\eps\eta\}$, and
$K_{cp}=\{\varphi\in K|\;\varphi$ is a convex combination of
continuous, normalized pure positive definite functions on
$G\,\}$. By \cite[Theorem C.5.5]{BHV}, $K_{cp}$ is weak* dense in
$\mathcal{P}_1(G)$. Since $K_{u,\eps\eta}$ is open in $K$, this
implies that $K_{cp}\cap K_{u,\eps\eta}$ is weak*-dense in
$\mathcal{P}_1(G)\cap K_{u,\eps\eta}$. By Raikov's Theorem, it is
also dense for the topology of uniform convergence on compact
subsets. Hence, since for all $\varphi\in K_{cp}\cap
K_{u,\eps\eta}$, $\varphi\ge 1-3\eps$ on $X$, the same holds for
all $\varphi\in \mathcal{P}_1(G)\cap K_{u,\eps\eta}$.\end{proof}

\begin{thm}
Let $G$ be a locally compact, $\sigma$-compact group. The
following are equivalent.
\begin{itemize}\item[(1)] $(G,X)$ has relative Property~(T).

\item[(2)] For every $\eps>0$, there exists a neighbourhood $V$ of
$1_G$ in $\hat{G}$ such that every $\pi\in V$ has a
$(X,\eps)$-invariant vector.

\item[(2')] For some $\eps<\sqrt{2}$, there exists a neighbourhood
$V$ of $1_G$ in $\hat{G}$ such that every $\pi\in V$ has a
$(X,\eps)$-invariant vector.\end{itemize}\label{thm:T_rep_irr}
\end{thm}
\begin{proof} (2)$\Rightarrow$(2') is trivial.

(2')$\Rightarrow$(1). By a result of Kakutani and Kodaira
\cite[Theorem~3.7]{Com}, there exists a compact normal subgroup
$K$ of $G$ such that $G/K$ is second countable. So we can suppose
that $G$ is second countable. Consider a unitary representation
$\pi$ of $G$ almost having invariant vectors. Arguing as in
\cite[proof of Lemme~1]{DK}, $\pi$ contains a nonzero
subrepresentation entirely supported by $V$. We conclude by Lemma
\ref{lem:desinteg_inv} below that $\pi$ has a $\eps'$-invariant
vector, where $\eps<\eps'<\sqrt{2}$. This proves that (3') of
Theorem \ref{thm:Del-Gui_relT} is satisfied.

(1)$\Rightarrow$(2). This is immediate from Condition (2) in
Theorem \ref{thm:Del-Gui_relT}.\end{proof}

\begin{rem}
The special case when $X$ is a subgroup is claimed without proof
in \cite[Chap. 1, 18.]{HV}.
\end{rem}

Let $G$ be a second countable, locally compact group, and
$X\subset G$. Let $(Z,\mu)$ be measured space, with $\mu(Z)>0$ and
$\mu$ $\sigma$-finite. Let $((\mathcal{H}_z)_{z\in Z},\Gamma)$ be
a measurable field of Hilbert spaces \cite[A 69]{Dix}, $\Gamma$
denoting a space of measurable vector fields. Let $(\pi_z)$ be a
field of unitary representations, meaning that $z\mapsto
\pi_z(g)x(z)$ is measurable, for every $x\in\Gamma$, $g\in G$.
Recall that, by definition, there exists a sequence $(x_n)$ in
$\Gamma$ such that, for every $z\in Z$, the family $(x_n(z))$ is
total in $\mathcal{H}_z$. Set $\pi=\int^\oplus\pi_zd\mu(z)$.

\begin{lem}
Fix $\eps>0$. Suppose that, for every $z$, $\pi_z$ has a
$(X,\eps)$-invariant vector. Then $\pi$ has a
$(X,\eps')$-invariant vector for every
$\eps'>\eps$.\label{lem:desinteg_inv}
\end{lem}

\begin{proof} Fix $0<\eta<1$. First note that replacing the family $(x_n)$
by the family of all its rational combinations if necessary, we
can suppose that, for every $z\in Z$ and every $v\in
\mathcal{H}_z$ of norm one, there exists $n$ such that
$\|v-x_n(z)\|\le\eta$. In particular, if $v$ is
$(X,\eps)$-invariant, then, for all $x\in X$,
$\|\pi_z(x)x_n(z)-x_n(z)\|\le \eps+2\eta$ and
$\|x_n(z)\|\ge\|1-\eta\|$, so that $x_n(z)$ is
$(X,(\eps+2\eta)/(1-\eta))$-invariant. Now define, for all $n\ge
0$,
$$A_n=\{z\in Z|\; 1-\eta\le\|x_n(z)\|\le 1+\eta,\; x_n(z)\text{ is
}(X,(\eps+2\eta)/(1-\eta))\text{-invariant}\}.$$

We have $\bigcup A_n=Z$ by the remark above. Using that $X$ is
separable, it is immediate that $A_n$ is measurable for every $n$.
Accordingly, there exists $n_0$ such that $\mu(A_{n_0})>0$. Using
that $\mu$ is $\sigma$-finite, there exists a measurable subset
$B\subset A_{n_0}$ such that $0<\mu(B)<\infty$. Define $\xi$ as
the
field $$z\mapsto\left\{%
\begin{array}{ll}
    x_{n_0}(z), & \hbox{$z\in B$;} \\
    0, & \hbox{otherwise.} \\
\end{array}%
\right.$$

Then it is clearly measurable, and
$$\|\xi\|^2=\int_{B}\|x_{n_0}(z)\|^2d\mu(z)\ge
(1-\eta)^2\mu(B),$$ and, for every $g\in X$,
$$\|\pi(g)\xi-\xi\|^2=\int_{B}\|\pi_z(g)x_{n_0}(z)-x_{n_0}(z)\|^2d\mu(z)$$ $$\le
((\eps+2\eta)^2(1+\eta)^2/(1-\eta)^2)\mu(B).$$ It follows that
$\xi\neq 0$ and is
$(X,(\eps+2\eta)(1+\eta)/(1-\eta)^2)$-invariant. Finally, for
every $\eps'>\eps$, we can choose $\eta$ sufficiently small so
that $(\eps+2\eta)(1+\eta)/(1-\eta)^2\le\eps'$.\end{proof}

\subsection{Some stability results}

We note for reference the following immediate but useful result:

\begin{prop}
Let $G$ be a locally compact group and $X_1,\dots,X_n$ be subsets.
Denote by $X_1\dots X_n$ the pointwise product $\{x_1\dots
x_n|\;(x_1,\dots,x_n)\in X_1\times\dots\times X_n\}$. Suppose
that, for every $i$, $(G,X_i)$ has relative Property~(T) [resp.
(FH)]. Then $(G,X_1\dots X_n)$ has relative Property~(T) [resp.
(FH)].\label{prop:relT_prod_subsets}
\end{prop}
\begin{proof} It suffices to prove the case when $n=2$, since then the
result follows by induction. For the case of Property (FH), this
follows from the inequality $\psi(gh)^{1/2}\le
\psi(g)^{1/2}+\psi(h)^{1/2}$ for every conditionally negative
definite function $\psi$. For the case of Property~(T), a similar
inequality holds since, if $\varphi$ is normalized positive
definite, then $1-|\varphi|^2$ is conditionally negative
definite.\end{proof}

\begin{exe}
1) If $G_1,\dots,G_n$ are locally compact groups, and $X_i\subset
G_i$, and if $(G_i,X_i)$ has relative Property~(T) [resp. (FH)]
for every $i$, then $(\prod G_i,\prod X_i)$ also has relative
Property~(T) [resp. (FH)].

2) Let $G$ be a group, $H_1$ a subgroup, and $H_2$ a subgroup of
finite index in $H_1$. If $(G,H_2)$ has relative Property~(T),
then so does $(G,H_1)$. It suffices to apply Proposition
\ref{prop:relT_prod_subsets} to $n=2$, $X_1=H_2$, and $X_2$ a
finite transversal of $H_1$ modulo $H_2$.

3) Fix $n\ge 3$, let $A$ a topologically finitely generated
locally compact commutative ring, and set $G=\SL_n(A)$. Denote by
$V_{n,m}$ the elements in $G$ that are products of $\le m$
elementary matrices. Then it follows from \cite[Corollary
3.5]{Sha99p} that $(G,V_{n,m})$ has relative Property~(T) for all
$m$. It is not known whether, for such $A$, there exists $m$ such
that $V_{n,m}=E_n(A)$, the subgroup generated by elementary
matrices; this seems to be an open question whenever $A$ is a
finitely generated commutative ring of Krull dimension $\ge 2$,
for instance $A=\Z[X]$, or $A=\mathbf{F}_p[X,Y]$.
\end{exe}

The following proposition is trivial.

\begin{prop}
Let $G$ be a locally compact group, $H$ a closed subgroup, and
$Y\subset X\subset H$ subsets. If $(H,X)$ has relative
Property~(T) [resp. (FH)], then so does $(G,Y)$.\qed
\end{prop}

\begin{prop}[Stability by extensions]
Let $G$ be a locally compact group, $N$ a closed normal subgroup,
and $X\subset G$ a subset. Denote by $p:G\to G/N$ the projection.

If $(G,N)$ and $(G/N,p(X))$ have relative Property~(T) [resp.
(FH)], then so does $(G,X)$.\label{prop:relT_extensions}
\end{prop}
\begin{proof} The assertion about relative Property (FH) is
immediate; that about relative Property~(T) is straightforward,
using the following fact\cite[Lemma B.1.1]{BHV}: for every compact
subset $K$ of $G/N$, there exists a compact subset $\tilde{K}$ of
$G$ such that $p(\tilde{K})=K$.\end{proof}

\subsection{Relative Property~(T) and compact generation}

It is well-known that a locally compact group with Property~(T) is
compactly generated. We generalize this result. The following
lemma is the easy part of such a generalization; we are going to
use it to prove something stronger.

\begin{lem}
Let $G$ be a locally compact group, and let $X\subset G$ be a
subset such that $(G,X)$ has relative Property~(T). Then $X$ is
contained in an open, compactly generated subgroup of
$G$.\label{lem:Trel_inc_cpt_eng}
\end{lem}
\begin{proof} For every open, compactly generated subgroup $\Omega$ of
$G$, let $\lambda_\Omega$ be the quasi-regular representation of
$G$ on $\ell^2(G/\Omega)$. Let $\delta_\Omega\in \ell^2(G/\Omega)$
be the Dirac function on $G/\Omega$. Let $\varphi_\Omega$ be the
corresponding coefficient. Then $\varphi_\Omega$ tends to 1,
uniformly on compact subsets, when $\Omega$ becomes big. By
relative Property~(T), the convergence is uniform on $X$, so that,
for some $\Omega$, we have $|1-\varphi_\Omega|<1$ on $X$. Since
$\varphi_\Omega$ has values in $\{0,1\}$, this implies that
$\delta_{\Omega}=\delta_{g\Omega}$ for all $g\in X$, that is,
$g\in\Omega$. Hence $X\subset\Omega$.\end{proof}

The following theorem shows that, in a certain sense, all the
information about relative Property~(T) lies within compactly
generated subgroups.

\begin{thm}
Let $G$ be a locally compact group, and $X\subset G$ a subset.
Then $(G,X)$ has relative Property~(T) if and only if there exists
an open, compactly generated subgroup $H$ such that $X\subset H$
and $(H,X)$ has relative
Property~(T).\label{thm:propT_and_cpt_gen}
\end{thm}

\begin{lem}
Fix $0<\eps<1$. Let $G$ be a locally compact group and $X$ a
subset. Then $(G,X)$ has relative Property~(T) if and only if, for
every net $(\varphi_i)$ of normalized, real-valued continuous
positive definite functions converging to 1 uniformly on compact
subsets, eventually $|\varphi_i|>\eps$ on
$X$.\label{lem:propT_coefs_quitte_zero}
\end{lem}
\begin{proof} The forward implication is trivial. Suppose that $(G,X)$ does
not have relative Property~(T). Then there exists a net
$(\varphi_i)$ of normalized, real-valued continuous positive
definite functions converging to 1 uniformly on compact subsets,
such that $\alpha=\sup_{i}\inf_{g\in X}\varphi_i(g)<1$. Then, for
some $n\in\N$, $\alpha^n<\eps$. Hence, $(\varphi_i^n)$ is a net of
normalized, continuous positive definite functions that converges
to 1 uniformly on compact subsets, but, for no $i$,
$|\varphi_i^n|>\eps$ on $X$.\end{proof}

\begin{proof}[Proof of Theorem \ref{thm:propT_and_cpt_gen}] By Lemma
\ref{lem:Trel_inc_cpt_eng}, there exists $\Omega\supset X$ an
open, compactly generated subgroup. Let $(K_i)$ be an increasing
net of open, relatively compact subsets, covering $G$, and denote
by $H_i$ the subgroup generated by $K_i$. We can suppose that
$\Omega\subset H_i$ for all $i$.

Suppose by contradiction that, for every $i$, $(H_i,X)$ does not
have Property~(T). Then, using Lemma
\ref{lem:propT_coefs_quitte_zero}, for all $i$ and all $n$, there
exists a normalized, continuous positive definite function
$\varphi_{i,n}$ on $H_i$, such that $\varphi_{i,n}\ge 1-2^{-n}$ on
$K_i$ and $\inf_{X}\varphi_{i,n}\le 1/2$. Since $H_i$ is open in
$G$, we can extend $\varphi_{i,n}$ to all of $G$, by sending the
complement of $H_i$ to $0$. It is clear that the net
$(\varphi_{i,n})$ tends to 1 uniformly on compact subsets of $G$,
but $\inf_X\varphi_{i,n}\le 1/2$. This contradicts that $(G,X)$
has relative Property~(T).\end{proof}

\subsection{$H$-metric}\label{subs:Hmetric}

First recall that a length function on a group $G$ is a function
$L:G\to\R_+$ satisfying the subadditivity condition $L(gh)\le
L(g)+L(h)$ for all $g,h$, and such that $L(1)=0$ and
$L(g)=L(g^{-1})$ for all $g\in G$. A length function defines a
(maybe non-separated) left-invariant metric on $G$ by setting
$d(g,h)=L(g^{-1}h)$.

Observe that if $L_1,L_2$ are two length functions, then so is
$L=\max(L_1,L_2)$. Indeed, we can suppose $L_1(gh)\ge L_2(gh)$.
Then $L(gh)=L_1(gh)\le L_1(g)+L_1(h)\le L(g)+L(h)$.

Also observe that a pointwise limit of length functions is a
length function. If follows that the upper bound of a family of
length functions, provided that it is everywhere finite, is a
length function.

\medskip

Now let $G$ be a locally compact, compactly generated group, and
$K$ a relatively compact, open generating subset. Define $\Psi_K$
as the upper bound of all (continuous, real-valued) conditionally
negative definite functions $\psi$ such that $\psi\le 1$ on $K$.
Recall that if $\psi$ is a real-valued conditionally negative
definite function, then $\psi^{1/2}$ is a length function. It
follows that $\Psi_K^{1/2}$ is a length function. It is easily
checked that it defines a separated metric on $G$, whose closed
balls are closed (for the initial topology). We call it the
$H$-\textit{metric}. It is easy to observe that if $K$ and $L$ are
two open, relatively compact generating subsets, then there exist
constants $\lambda,\lambda'>0$ such that $\lambda\Psi_K\le
\Psi_L\le \lambda'\Psi_K$. Accordingly, the identity map defines a
bi-Lipschitz map between these two metrics, and the choice of $K$
is not essential at all.

\begin{prop}
Let $G$ be a locally compact, compactly generated group, and $X$ a
subset. Then $(G,X)$ has relative Property~(T) if and only if $X$
is bounded for the $H$-metric.
\end{prop}
\begin{proof} First recall that, since $G$ is $\sigma$-compact, relative
Property~(T) and relative Property (FH) are equivalent by Theorem
\ref{thm:Del-Gui_relT}.

If $X$ is bounded for the $H$-metric, and $\psi$ is a (continuous,
real-valued) conditionally negative definite function on $X$,
then, for some constant $\alpha>0$, $\alpha\psi\le 1$ on $K$. So
$\psi\le\alpha^{-1}\Psi_K$, which is bounded on $X$, and thus
$(G,X)$ has relative Property (FH).

Conversely, suppose that $X$ is not bounded for the $H$-metric.
Then there exist a sequence of (continuous, real-valued)
conditionally negative definite functions $\psi_n$, bounded by 1
on $K$, and a sequence $x_n$ of $X$ such that $\psi_n(x_n)\ge
4^n$. Set $\psi=\sum 2^{-n}\psi_n$. Since the convergence is
uniform on compact subsets, $\psi$ is a well-defined continuous
conditionally negative definite function on $G$, and $\psi(x_n)\ge
2^n$, so that $\psi$ is not bounded on $X$, and $(G,X)$ does not
have relative Property (FH).\end{proof}

\begin{cor}
Let $G$ be a locally compact, compactly generated group.

\begin{itemize}\item[1)] $G$ has Property~(T) is and only if it is bounded for the
$H$-metric.

\item[2)] $G$ has no non-compact closed subset with relative
Property~(T) if and only if $G$ is proper for the $H$-metric (that
is, the balls for the $H$-metric are compact for the initial
topology).\end{itemize}
\end{cor}

It is maybe interesting comparing the $H$-metric with the word
metric (relative to any compact generating set). A general result
in this direction has been obtained by V.~Lafforgue \cite{Laf}: if
$G$ does not have Property~(T), if $L_K$ denotes the word length
with respect to the compact generating set $K$, then there exists
a conditionally negative definite function $\psi$ on $G$
satisfying $\psi|_K\le 1$ and
$$\sup\{\psi(x)^{1/2}|\,x\in G\hbox{ and }L_K(x)\le n\}
\ge\frac{\sqrt{n}}{2}-2.$$ In
particular, if $H_K=\Psi_K^{1/2}$ is the length in the $H$-metric,
then, for every $0<C<1/2$, there exists a sequence $(x_n)$ in $G$
such that, for all $n$, $L_K(x_n)\le n$ and $H_K(x_n)\ge
C\sqrt{n}$.

On the other hand, R. Tessera \cite{Tes} has proved that, for a
polycyclic group $\Gamma$, the $H$-metric is equivalent to the
word metric, although it is known that $\Gamma$ has no
conditionally negative definite function with quadratic growth
unless $\Gamma$ is virtually abelian~\cite{CTV}.

\section{Relative Property~(T) in connected Lie groups and $p$-adic algebraic
groups}\label{Sec_relT_lie_alg}

\subsection{Preliminaries}

Given a locally compact group $G$, we can naturally raise the
problem of determining for which subsets $X$ the pair $(G,X)$ has
relative Property~(T). Here is a favourable case, where the
problem is completely solved.

\begin{lem}
Let $G$ be a locally compact group, and $N$ a normal subgroup such
that $(G,N)$ has relative Property~(T) and $G/N$ is Haagerup. Let
$X$ be any subset of $G$. Then $(G,X)$ has relative Property~(T)
if and only if the image of $X$ in $G/N$ is relatively
compact.\label{lem:favour_Trel}
\end{lem}
\begin{proof} The condition is clearly necessary, since relative
Property~(T) is inherited by images.

Conversely, if the image of $X$ in $G/N$ is relatively compact,
there exists a compact subset $K$ of $G$ such that $X$ is
contained in $KN=\{kn|\;(k,n)\in K\times N\}$. Let $\psi$ be a
continuous, conditionally negative definite function on $G$. Then
$\psi$ is bounded on $N$ and on $K$, hence on $KN$, hence on $X$.
This proves that $(G,X)$ has relative Property (FH). In view of
Theorem \ref{thm:Del-Gui_relT}, this is sufficient if $G$ is
$\sigma$-compact. Actually, we can reduce to this case: indeed, by
Theorem \ref{thm:propT_and_cpt_gen}, there exists an open,
compactly generated subgroup $H$ of $G$, that contains $N$ and can
be supposed to contain $K$, such that $(H,N)$ has relative
Property~(T).\end{proof}

Recall the key result, due to Shalom \cite[Theorem~5.5]{Sha99t}
(see also \cite[\S 1.4]{BHV}).

\begin{prop}
Let $G$ be a locally compact group and $N$ a closed normal abelian
subgroup. Assume that the only mean on the Borel subsets of the
Pontryagin dual $\widehat{N}=\textnormal{Hom}(N,\R/\Z)$, invariant
under the action of $G$ by conjugation, is the Dirac measure at
zero. Then the pair $(G,N)$ has relative
Property~(T).\label{dirac}
\end{prop}

This result allows to prove relative Property~(T) for certain
normal abelian subgroups. Since we also deal with nilpotent
subgroups, we use the following proposition, which generalizes
\cite[Proposition 4.1.4]{CCJJV}.

\begin{prop}
Let $G$ be a locally compact, $\sigma$-compact group, $N$ a closed
subgroup, and let $Z$ be a closed, central subgroup of $G$
contained in $\overline{[N,N]}$. Suppose that every morphism of
$N$ into a compact Lie group has an abelian image.

Suppose that the pair $(G/Z,N/Z)$ has Property~(T). Then $(G,N)$
has Property~(T).\label{prop:Trel_ext_centrale}
\end{prop}
\begin{proof} It suffices to show that $(G,Z)$ has relative Property~(T).
Indeed, since the pairs $(G,Z)$ and $(G/Z,N/Z)$ have relative
Property~(T), it then follows by Proposition
\ref{prop:relT_extensions} that $(G,N)$ has relative Property~(T).

We use an argument similar to the proof of \cite[Lemma 1.6]{Wang}.
To show that $(G,Z)$ has relative Property~(T), we use the
characterization by nets of irreducible unitary representations
(see Theorem \ref{thm:T_rep_irr}). Let $\pi_i$ be a net of
irreducible unitary representations of $G$ converging to the
trivial representation: we must show that eventually $\pi_i$
factors through $Z$. Let $\overline{\pi}_i$ be the contragredient
representation of $\pi_i$. Then $\pi_i\otimes\overline{\pi_i}$
converges to the trivial representation. By irreducibility,
$\pi_i$ is scalar in restriction to $Z$, hence
$\pi_i\otimes\overline{\pi_i}$ is trivial on $Z$, so factors
through $G/Z$. Since $(G/Z,N/Z)$ has Property~(T), the restriction
to $N$ of $\pi_i\otimes\overline{\pi_i}$ eventually contains the
trivial representation. By a standard argument \cite[Appendix
1]{BHV}, this means that $\pi_i|_N$ eventually contains a
finite-dimensional subrepresentation~$\rho_i$.

Remark that $\overline{\rho_i(N)}$ is a compact Lie group; so it
is, by assumption, abelian. This means that $[N,N]$ acts
trivially; hence $Z$ does so as well: $\rho_i$ is trivial on $Z$.
Hence, for large $i$, $\pi_i$ has nonzero $Z$-invariants vectors;
by irreducibility, $\pi_i$ is trivial on $Z$. Accordingly $(G,Z)$
has Property~(T).\end{proof}

We shall use the following well-known result of Furstenberg
\cite{FUR}.

\begin{thm}[Furstenberg]
Let $\K$ be a local field, $V$ a finite dimensional $K$-vector
space. Let $G\subset \textnormal{PGL}(V)$ be a Zariski connected
(but not necessarily Zariski closed) subgroup, whose closure is
not compact. Suppose that $G$ preserves a probability measure
$\mu$ on the projective space $\mathbf{P}(V)$. Then there exists a
$G$-invariant proper projective subspace $W\subsetneq
\mathbf{P}(V)$ such that $\mu(W)=1$.\label{furcon}
\end{thm}

\begin{rem}
Observe that a subgroup of $\textnormal{PGL}(V)$ preserves an
invariant mean on $\mathbf{P}(V)$ if and only if it preserves a
probability: indeed, a mean gives rise to a normalized positive
linear form on $L^\infty(\mathbf{P}(V))$, and restricts to a
normalized positive linear form on $C(\mathbf{P}(V))$, defining a
probability.\label{rem:cpt_mean_proba}
\end{rem}

We say that a topological group $G$ is discompact\footnote{This is
often called ``minimally almost periodic", but we prefer the
terminology ``discompact", introduced in \cite{Sha99t}.} if there
is no nontrivial morphism of $G$ to a compact group.

\begin{rem}
If $G$ is a discompact locally compact group, then it has trivial
abelianization. Indeed, it follows that its abelianization is also
discompact, so has trivial Pontryagin dual, so it trivial by
Pontryagin duality.\label{rem:discompact_abelianization}
\end{rem}

The following corollary is an easy consequence of Theorem
\ref{furcon} (see \cite[Corollary~2.39]{CorThe} for a proof).

\begin{cor}
Let $G$ be a discompact locally compact group. Let $V$ be a
finite-dimensional vector space over $\K$, and let $G\to
\textnormal{GL}(V)$ be any continuous representation. Then $G$
preserves a probability on $\mathbf{P}(V)$ if and only if $G$ has
a nonzero fixed point on $V$.\qed\label{cor:no_proba_preserved}
\end{cor}

\begin{exe}
(1) Let $G$ be a simply connected, simple group over $\K$, of
positive $\K$-rank. Then $G(\K)$ is discompact. Indeed, $G(\K)$ is
generated by elements whose conjugacy classes contain 1 in their
closure: this follows from the following observations: $G(\K)$ is
simple \cite[Chap.~I, Theorem~1.5.6 and
Theorem~2.3.1(a)]{Margulis}, and there exists a subgroup of $G$
isomorphic to either $\SL_2(\K)$ or $\textnormal{PSL}_2(\K)$
\cite[Chap.~I, Proposition 1.6.3]{Margulis}. Accordingly, every
morphism of $G(\K)$ into a compact group has trivial image.

(2) Let $G$ be a connected, noncompact, simple Lie group. Then $G$
is discompact. Indeed, such a group is generated by connected
subgroups locally isomorphic to $\SL_2(\R)$, hence is generated by
elements whose conjugacy class contains 1.\label{exe:discompact}
\end{exe}

\begin{prop}
Let $V$ be a finite-dimensional space over a local field $\K$. Let
$G$ be any locally compact group, and $\rho:G\to
\textnormal{GL}(V)$ a continuous representation. Then $(G\ltimes
V,V)$ has relative Property~(T) if and only if $G$ preserves no
probability on $\mathbf{P}(V^*)$. In particular, if $G$ is
discompact, then $(G\ltimes V,V)$ has relative Property~(T) if and
only if $G$ fixes no point in $V^*$. \label{prop:propT_proba_proj}
\end{prop}

\begin{proof} Suppose, by contradiction that $G$ preserves a probability on
$\mathbf{P}(V^*)$ and $(G\ltimes V,V)$ has relative Property~(T).
By Theorem \ref{furcon}, the finite index subgroup $G_0$ of $G$
(its unit component in the inverse image of the Zariski topology
from $\textnormal{GL}(V)$) preserves a nonzero subspace $W\subset
V^*$, such that the image of the morphism $G_0\to
\textnormal{PGL}(W)$ has compact closure. Since $W$ is a subspace
of $V^*$, $W^*$ is a quotient of $V$. By Corollary 4.1(2) in
\cite{Jol}, $(G_0\ltimes V,V)$ has relative Property~(T), and so
has $(G_0\ltimes W^*,W^*)$. This implies that
$(\overline{\rho(G_0)}\ltimes W^*,W^*)$ also has relative
Property~(T). But $\overline{\rho(G_0)}\ltimes W^*$ is amenable,
so that $W^*$ is compact, and this is a contradiction.

The converse is due to M. Burger \cite[Proposition~7]{Bur}. It is
obtained by combining Proposition \ref{dirac} and Remark
\ref{rem:cpt_mean_proba}. The second assertion follows from
Corollary \ref{cor:no_proba_preserved}.\end{proof}

\begin{rem}
It is worth noting that, in Proposition
\ref{prop:propT_proba_proj}, and in view of Corollary
\ref{cor:no_proba_preserved}, relative Property~(T) for $(G\ltimes
V,V)$ only depends on the closure (for the ordinary topology) of
the image of $G$ in $\textnormal{PGL}(V)$.
\end{rem}

\subsection{Relative Property~(T) in
algebraic groups over local fields of characteristic
zero}\label{subs:relT_alg}

We denote by $\K$ a local field of characteristic zero. Here is
the main lemma of this paragraph.

\begin{lem}
Let $G$ be a linear algebraic $\K$-group, which decomposes as
$S\ltimes R$, where $S$ is semisimple and $\K$-isotropic, and $R$
is unipotent.

Suppose that $[S,R]=R$. Then $(G(\K),R(\K))$ has relative
Property~(T).\label{lem:relT_alg}
\end{lem}
%Replacing $G$ by $S_{mn}\ltimes R_u$, which is a characteristic
%subgroup of $G$, we can suppose that $S=S_{mn}$.
\begin{proof} Replacing $S$ by its universal cover if necessary, we can
suppose that $S$ is simply connected. We then argue by induction
on the dimension of $R$. If the dimension is zero, there is
nothing to prove; suppose $R\neq 1$. Let $Z$ be the last nonzero
term of its descending central series.

First case: $Z$ is central in $G$. The hypothesis $[S,R]=R$
implies that $R$ is not abelian. Hence $Z\subset [R,R]$, so that
$Z(\K)\subset[R,R](\K)$. By \cite[Lemma 13.2]{BT}
$[R,R](\K)=[R(\K),R(\K)]$, so that $Z(\K)\subset[R(\K),R(\K)]$. We
must check that the hypotheses of Proposition
\ref{prop:Trel_ext_centrale} are fulfilled. Let $W$ be a compact
Lie group, and $R(\K)\to W$ a morphism with dense image: we must
show that $W$ is abelian. Since $R(\K)$ is solvable, the connected
component $W_0$ is abelian. Moreover, $R(\K)$ is divisible, so
$W/W_0$ is also divisible; this implies $W=W_0$. Accordingly, by
Proposition \ref{prop:Trel_ext_centrale}, since
$(G(\K)/Z(\K),R(\K)/Z(\K))$ has relative Property~(T) by induction
hypothesis, it follows that $(G(\K),R(\K))$ has relative
Property~(T).

Second case: $Z$ is not central in $G$. Set $N=[S,Z]$. Then
$[S,N]=N$. By Proposition \ref{prop:propT_proba_proj} and in view
of Example \ref{exe:discompact}(1), $(G(\K),N(\K))$ has relative
Property~(T). By the induction assumption,
$((G/N)(\K),(R_u/N)(\K))$, which coincides with
$(G(\K)/N(\K),R_u(\K)/N(\K))$, has relative Property~(T). Hence
$(G(\K),R_u(\K))$ has relative Property~(T).\end{proof}

Let $G$ be a linear algebraic group over $\K$. We denote by $R_u$
its unipotent radical, and $L$ a reductive Levi factor (so that
$G_0=L\ltimes R_u$). We decompose $L$ as an almost product
$L_mL_{nm}$, where $L_m$ (resp. $L_{nm}$) includes the centre of
$L$, and the simple factors of rank zero (resp. includes the
simple factors of positive rank)\footnote{$(n)m$ stands for
(non-)amenable.}.

Let $R$ be the radical of $G$, let $S$ be a Levi factor, and
decompose it as $S_{c}S_{nc}$, where $S_c$ (resp. $S_{nc}$ is the
sum of all factors of rank 0 (resp. of positive rank).

If $\g$ is a Lie algebra and $\mk{h}_1$, $\mk{h}_2$ are two
subspaces, we denote by $[\mk{h}_1,\mk{h}_2]$ (resp.
$[\mk{h}_1,\mk{h}_2]_v$) the Lie algebra (resp. the subspace)
generated by the elements of the form $[h_1,h_2]$,
$(h_1,h_2)\in\mk{h}_1\times\mk{h}_2$.

\begin{lem}
For every Levi factors $L,S$ of respectively $R_u$ and $R$, we
have $[L_{nm},R_u]=[S_{nc},R]$, and this is a $\K$-characteristic
subgroup of $G$.\label{lem:ScnR=LnmRu}
\end{lem}
\begin{proof} We can work within the Lie algebra. We first justify that
$[\mk{l}_{nm},\mk{r}_u]$ is an ideal: indeed,
$$[\mk{l},[\mk{l}_{nm},\mk{r}_u]]\subset[[\mk{l},\mk{l}_{nm}],\mk{r}_u]+[\mk{l}_{nm},[\mk{l},\mk{r}_u]]
\subset [\mk{l}_{nm},\mk{r}_u]$$ since
$[\mk{l},\mk{l}_{nm}]\subset\mk{l}_{nm}$ and
$[\mk{l},\mk{r}_u]\subset\mk{r}_u$. On the other hand,
\begin{eqnarray*}
  [\mk{r}_u,[\mk{l}_{nm},\mk{r}_u]] &=& [\mk{r}_u,[\mk{l}_{nm},[\mk{l}_{nm},\mk{r}_u]]] \\
   &\subset &
   [\mk{l}_{nm},[\mk{r}_u,[\mk{l}_{nm},\mk{r}_u]]]\;+
   \;[[\mk{l}_{nm},\mk{r}_u],[\mk{l}_{nm},\mk{r}_u]]\\
   &\subset & [\mk{l}_{nm},\mk{r}_u].
\end{eqnarray*}

It follows that $[L_{nm},R_u]$ is a normal subgroup of $G$. By
\cite[(5.1)]{BS}, the $\K$-conjugacy class of $L$ does not depend
of the choice of $L$. So the same thing holds for $L_{nm}$ (which
is $\K$-characteristic in $L$). Accordingly, $[L_{nm},R_u]$ is a
$\K$-characteristic subgroup of $G$.

Now, since $S_{nc}$ is a reductive $\K$-subgroup of $G$, again
using \cite[(5.1)]{BS}, up to $\K$-conjugate if necessary, we can
suppose that $S_{nc}\subset L$, so that finally $S_{nc}=L_{nm}$,
and $R=L_r\ltimes R_u$, where $L_r$ is the unit component of
centre of $L$. Since $[L_{nc},L_r]=1$, we obtain
$[S_{nc},R]=[L_{nc},R]=[L_{nc},R_u]$.\end{proof}

Let $S_{nh}$ be the sum of all simple factors $H$ of $S_{nc}$ such
that $H(\K)$ is not Haagerup (equivalently: has Property~(T)):
these are factors of rank $\ge 2$, and also, when $\K=\R$, factors
locally isomorphic to $\textnormal{Sp}(n,1)$ or $F_{4(-20)}$.

\begin{defn}Define $R_T$ as the $\K$-subgroup
$S_{nh}[S_{nc},R]$ of $G$.\end{defn}

\begin{thm}
$R_T$ is a $\K$-characteristic subgroup of $G$, the quotient group
$G(\K)/R_T(\K)$ is Haagerup, and \allowbreak$(G(\K),R_T(\K))$ has
relative Property~(T).\label{thm:RT_alg}
\end{thm}
\begin{proof} It follows from Lemma \ref{lem:ScnR=LnmRu} that $[S_{nc},R]$
is unipotent. Consider the $\K$-subgroup $W=S_{nc}[S_{nc},R]$ of
$G$. Applying Lemma \ref{lem:relT_alg} to $W$, we obtain that
$(G(\K),[S_{nc},R](\K))$ has relative Property~(T). Since
$(G(\K),S_{nh}(\K))$ also has relative Property~(T) and
$[S_{nc},R](\K)$ is a normal subgroup, we obtain that
$(G(\K),R_T(\K))$ has relative Property~(T) by Proposition
\ref{prop:relT_prod_subsets}.

To show that $R_T$ is a $\K$-characteristic subgroup, we can work
modulo the subgroup $[S_{nc},R]$ which is $\K$-characteristic by
Lemma \ref{lem:ScnR=LnmRu}. But, in $G/[S_{nc},R]$, $S_{nc}$ is a
direct factor and can be characterized as the biggest normal
subgroup that is connected, semisimple, and $\K$-isotropic; and
$S_{nh}$ is $\K$-characteristic in $S_{nc}$. It follows that $R_T$
is $\K$-characteristic.

Finally, $H=G/R_T$ is almost the direct product of a semisimple
group $H_s$ such that $H_s(\K)$ is Haagerup, and its amenable
radical $H_m$, such that $H_m(\K)$ is amenable, hence Haagerup. So
$H(\K)$ is Haagerup, and contains $G(\K)/R_T(\K)$ as a closed
subgroup.\end{proof}

So we are in position to apply Lemma \ref{lem:favour_Trel}.

\begin{cor}
Let $X$ be a subset of $G(\K)$. Then $(G(\K),X)$ has relative
Property~(T) if and only if the image of $X$ is $G(\K)/R_T(\K)$ is
relatively compact.\qed
\end{cor}

We retrieve a result of Wang (his statement is slightly different
but equivalent to this one).

\begin{cor}[Wang]
$G(\K)$ has Property~(T) if and only if $S_{nh}[S_{nc},R_u](\K)$
is cocompact in $G(\K)$.\qed
\end{cor}

\begin{cor}[\cite{CorJLT}]
$G(\K)$ is Haagerup if and only if $S_{nh}=[S_{nc},R]=1$.\qed
\end{cor}

\subsection{Relative Property~(T) in connected Lie groups}\label{subs:relT_lie}

Let $G$ be a Lie group (connected, even if it is straightforward
to generalize what follows to a Lie group with finitely many
components), $R$ its radical, $S$ a Levi factor (not necessarily
closed), decomposed as $S_cS_{nc}$ by separating compact and
noncompact factors. Let $S_{nh}$ be the sum of all simple factors
of $S_{nc}$ that have Property~(T). Set
$R_T=\overline{S_{nh}[S_{nc},R]}$.

\begin{thm}
$R_T$ is a characteristic subgroup of $G$, $G/R_T$ is Haagerup,
and $(G,R_T)$ has relative Property~(T).\label{thm:Trel_Lie}
\end{thm}

\begin{proof} The first statement can be proved in the same lines as in the
algebraic case.

It is immediate that $G/R_T$ is locally isomorphic to a direct
product $M\times S$ where $M$ is amenable and $S$ is semisimple
with all simple factors locally isomorphic to $\SO(n,1)$ or
$\SU(n,1)$. By \cite[Chap. 4]{CCJJV}, $G/R_T$ is Haagerup.

Finally, let us show that $(G,R_T)$ has relative Property~(T).
First, note that we can reduce to the case when $G$ is simply
connected. Indeed, let $p:\tilde{G}\to G$ be the universal
covering. Then $p(\tilde{H})=H$, for $H=R,S_{nc},S_{nh}$, where
$\tilde{H}$ is the analytic subgroup of $\tilde{G}$ that lies over
$H$. If the simply connected case is done, then
$(\tilde{G},\widetilde{S_{nh}}[\widetilde{S_{nc}},\tilde{R}])$ has
relative Property~(T). It follows that
$(G,p(\widetilde{S_{nh}}[\widetilde{S_{nc}},\tilde{R}]))$ also has
relative Property~(T), and the closure of
$p(\widetilde{S_{nh}}[\widetilde{S_{nc}},\tilde{R}])$ is equal to
$R_T$.

Now suppose that $G$ is simply connected. Then the subgroup
$S_{nc}[S_{nc},R]$ is closed and isomorphic to $S_{nc}\ltimes
[S_{nc},R]$. Arguing as in the proof of Lemma \ref{lem:relT_alg}
(using Example \ref{exe:discompact}(2) instead of (1)),
$(S_{nc}\ltimes [S_{nc},R],[S_{nc},R])$ has relative Property~(T).
Hence $(G,[S_{nc},R])$ also has relative Property~(T), and, as in
the proof of Theorem \ref{thm:RT_alg}, it implies that
$(G,S_{nh}[S_{nc},R])$ has relative Property~(T).\end{proof}

So we are again in position to apply Lemma \ref{lem:favour_Trel}.

\begin{cor}
Let $X$ be a subset of $G$. Then $(G,X)$ has relative Property~(T)
if and only if the image of $X$ in $G/R_T$ is relatively compact.
\end{cor}

We also retrieve a result of Wang in the case of connected Lie
groups.

\begin{cor}[Wang]
The connected Lie group $G$ has Property~(T) if and only if
$\overline{S_{nh}[S_{nc},R_u]}$ is cocompact in $G$.\qed
\end{cor}

\begin{cor}[{\cite[Chap.~4]{CCJJV}}]
The connected Lie group $G$ is Haagerup if and only
if\footnote{Note that we used, in the proof of Theorem
\ref{thm:Trel_Lie}, one implication from \cite[Chap. 4]{CCJJV},
namely, $S_{nh}=[S_{nc},R]=1$ implies $G$ Haagerup. This result is
easy when $S_{nc}$ has finite centre, but, otherwise, is much more
involved. Accordingly, only the reverse implication can be
considered as a corollary of the present work.}
$S_{nh}=[S_{nc},R]=1$.\label{cor:Lie_haag}
\end{cor}

\begin{proof}[Proof of Corollary \ref{cor:Lie_haag}] The hypothesis implies that
$S_{nh}$ and $W=[S_{nc},R]$ are both relatively compact. So
$S_{nh}=1$. Now, since $[S_{nc},[S_{nc},R]]=[S_{nc},R]$, we have
$[S_{nc},\overline{W}]=\overline{W}$. But, since $W$ is a compact,
connected, and solvable Lie group, it is a torus; since $S_{nc}$
is connected, its action on $\overline{W}$ is necessarily trivial,
so that $W\subset[S_{nc},\overline{W}]=1$.\end{proof}

\begin{rem}
If $G$ is a connected Lie group without the Haagerup Property, the
existence of a noncompact closed subgroup with relative
Property~(T) was proved in \cite{CCJJV}, and later established by
another method in \cite{CorJLT}, where the result was generalized
to linear algebraic groups over local fields of characteristic
zero. However, in both cases, the subgroup constructed is not
necessarily normal, while $R_T$ is.
\end{rem}

\begin{rem}
In this remark, given a locally compact group $G$, we say that a
closed, normal subgroup $N$ is a T-radical if $G/N$ is Haagerup
and $(G,N)$ has relative Property~(T).

It is natural to ask about the uniqueness of T-radicals when they
exist. Observe that if $N,N'$ are T-radicals, then the image of
$N$ in $G/N'$ is relatively compact, and vice versa. In
particular, if $G$ is discrete, then all T-radicals are
commensurable. This is no longer the case if $G$ is not discrete,
for instance, set $G=\textnormal{SL}(2,\Z)\ltimes\R^2$. Then the
subgroups $a\Z^2$, for $a\neq 0$, are all T-radicals, although two
of them may have trivial intersection.

In this example, $G$ has no minimal T-radical. This is also the
case in $\SL(2,\Z)\ltimes\Z^2$. On the other hand, let $G$ be a
finitely generated solvable group with infinite locally finite
centre. Then, every finite subgroup of the centre is a T-radical,
but $G$ has no infinite T-radical, so has no maximal T-radical. An
example of such a group $G$ is the group of matrices of the form
$\begin{pmatrix}
  1 & a & b \\
  0 & u^n & c \\
  0 & 0 & 1 \
\end{pmatrix}$, for $a,b,c\in\mathbf{Z}[u,u^{-1}]$, $n\in\Z$.

However, if $G$ is a connected Lie group, it can be shown that $G$
has a minimal and a maximal T-radical. The minimal one is $R_T$,
as defined above: indeed, if $H$ is a quotient of $G$ with the
Haagerup Property, then $S_{nh}$ and $[S_{nc},R]$ are necessarily
contained in the kernel. The maximal one is found by taking the
preimage of the maximal normal compact subgroup of $G/R_T$; it can
immediately be generalized to any connected locally compact group.
\end{rem}

\begin{rem}
Following Shalom \cite{Sha99t}, if $G$ is a topological group and
$H$ is a subgroup, we say that $(G,H)$ has {\em strong} relative
Property~(T) if there exists a Kazhdan pair $(K,\eps)$ for the
pair $(G,H)$ with $K$ finite (and $\eps>0$). More precisely, this
means that every unitary representation with a
$(K,\eps)$-invariant vector has a $H$-invariant vector. In this
context, it is natural to equip $\hat{G}$ with the topology
inherited from $\widehat{G_d}$, the unitary dual of $G_d$, where
$G_d$ denotes $G$ with the discrete topology. As for the case of
relative Property~(T), it can be checked that $(G,H)$ has strong
relative Property~(T) if and only if, for every net $\pi_i$ in
$\widehat{G}$ converging to 1 in $\widehat{G_d}$, eventually
$\pi_i$ has a $H$-invariant vector. Then it is straightforward
from the proof that Proposition \ref{prop:Trel_ext_centrale}
remains true for strong relative Property~(T). On the other hand,
Proposition \ref{dirac} is actually true with strong Property~(T)
\cite[Theorem~5.5]{Sha99t}. It then follows from the proofs above
that, if $G$ is a connected Lie group, then $(G,R_T)$ has strong
relative Property~(T), and similarly for algebraic groups over
local fields of characteristic zero.
\end{rem}

%%%%%%%%%%%%%%%%%%%%%%%%%%%%%%%%%%%%%%%%%%%%%%%%%%%%%%%%%%%%%%%%%%%%%%%%%%%%%%%

\section{Framework for irreducible lattices: resolutions}\label{Sec_resol}

In this section, we make a systematic study of ideas relying on
work of Lubotzky and Zimmer \cite{LZ}, and later apparent in
\cite[Chap. III, 6.]{Margulis} and \cite{BL}.

Given a locally compact group $G$, when can we say that we have a
good quantification of Kazhdan's Property~(T)? Lemma
\ref{lem:favour_Trel} provides a satisfactory answer whenever $G$
has a normal subgroup $N$ such that $G/N$ is Haagerup and $(G,N)$
has relative Property~(T). We have seen in Section
\ref{Sec_relT_lie_alg} that this is satisfied in a large class of
groups. However, this is not inherited by lattices. A typical
example is the case of an irreducible lattice $\Gamma$ in a
product of noncompact simple connected Lie groups $G\times H$,
where $G$ has Property~(T) and $H$ is Haagerup. In such an
example, although $\Gamma\cap G=\{1\}$, $G$ can be thought as a
``ghost" normal sugroup of $\Gamma$, and is the ``kernel" of the
projection $\Gamma\to H$. Relative Property~(T) for the pair
$(G\times H,G)$ can be restated by saying that the projection
$G\times H\to H$ is a ``resolution". By a theorem essentially due
to Margulis (Theorem \ref{thm:MBL}), this notion is inherited by
lattices, so that, in this case, the projection $\Gamma\to H$ is a
resolution.

Before giving rigorous definitions, we need some elementary
preliminaries.

\subsection{$Q$-points}

We recall that an action $\alpha$ by isometries of a topological
group $G$ on a metric space $X$ is continuous if the function
$g\mapsto\alpha(g)x$ is continuous for every $x\in X$. All the
functions and actions here are supposed continuous.

Let $f:G\to Q$ be a morphism between topological groups, with
dense image. Recall that, for any Hausdorff topological space $X$,
a function $u:G\to X$ factors through $Q$ if and only if, for
every net $(g_i)$ in $G$ such that $f(g_i)$ converges in $Q$, the
net $u(g_i)$ converges in $X$; note that the factorization $Q\to
X$ is unique.

\begin{defn}Let $f:G\to Q$ be a morphism between topological groups, with
dense image. Let $\alpha$ be an action of $G$ by isometries on a
metric space $X$. We call $x\in X$ a $Q$-point if the orbital map
$g\mapsto\alpha(g)x$ factors through $Q$.\end{defn}

\begin{prop}
The set $X^Q$ of $Q$-points in $X$ is $G$-invariant, and the
action $\alpha^Q$ of $G$ on $X^Q$ factors through $Q$. If moreover
$X$ is a complete metric space, then $X^Q$ is closed in
$X$.\label{prop:H^Q_closed}
\end{prop}
\begin{proof} The first assertion is immediate; let us assume that
$X$ is complete and let us show that $X^Q$ is closed. Let $(y_n)$
be a sequence in $X^Q$, converging to a point $y\in X$. Write
$\alpha(g)y_n=w_n(f(g))$, where $w_n$ is a continuous function:
$Q\to X$. If $m,n\in \N$,
$d(w_m(f(g)),w_n(f(g)))=d(\alpha(g)y_m,\alpha(g)y_n)=d(y_m,y_n)$.
It follows that $\sup_{q\in f(G)}d(w_m(q),w_n(q))\to 0$ when
$m,n\to\infty$. On the other hand, since $f(G)$ is dense in $Q$,
$\sup_{q\in f(G)}d(w_m(q),w_n(q))=\sup_{q\in Q}d(w_m(q),w_n(q))$.
Accordingly, $(w_n)$ is a Cauchy sequence for the topology of
uniform convergence on $Q$. Since $X$ is complete, this implies
that $(w_n)$ converges to a continuous function $w:Q\to X$.
Clearly, for all $g\in G$, $\alpha(g)y=w(f(g))$, so that $y\in
X^Q$.\end{proof}

\begin{prop}
Suppose that $X$ is a complete metric space. Given $x\in X$, the
following are equivalent:

\begin{itemize}\item[1)] $x\in X^Q$.

\item[2)] The mapping $g\mapsto d(x,\alpha(g)x)$ factors through
$Q$.

\item[3)] For every net $(g_i)$ in $G$ such that $f(g_i)\to 1$,
$d(x,\alpha(g_i)x)\to 0$.\end{itemize}\label{prop:H^Q_fct}
\end{prop}
\begin{proof} 1)$\Rightarrow$2)$\Rightarrow$3) is immediate.

Suppose 3). Let $(g_i)$ be a net in $G$ such that $f(g_i)$
converges in $Q$. Then $(f(g_i^{-1}g_j))$ converges to 1 when
$i,j\to\infty$, so that
$$d(\alpha(g_i)x,\alpha(g_j)x)=d(y,\alpha(g_i^{-1}g_j)x)\to 0,$$
i.e. $(\alpha(g_i)x)$ is Cauchy. Hence, it converges since $X$ is
complete. This means that $g\mapsto\alpha(g)x$ factors through
$Q$, i.e. $x\in X^Q$.\end{proof}

Recall that CAT(0) metric spaces are a generalization of simply
connected Riemannian manifolds with non-positive curvature; see
\cite{BH} for a definition.

\begin{prop}~
\begin{itemize}\item[1)] Suppose that $X$ is a complete CAT(0) metric space. Then $X^Q$ is
a closed, totally geodesic subspace.

\item[2)] If $X=\mathscr{H}$ is the Hilbert space of a unitary
representation $\pi$ of $G$, then $\mathscr{H}^Q$ is a closed
subspace, defining a subrepresentation $\pi^Q$ of $\pi$ (we refer
to elements in $\mathscr{H}^Q$ as $Q$-vectors rather that
$Q$-points). For every $\xi\in\mathscr{H}$, $\xi\in\mathscr{H}^Q$
if and only if the corresponding coefficient
$g\mapsto\langle\xi,\pi(g)\xi\rangle$ factors through $Q$.

\item[3)] If $X=\mathscr{H}$ is an affine Hilbert space, then
$\mathscr{H}^Q$ is a closed affine subspace (possibly empty). For
every $v\in\mathscr{H}$, $v\in\mathscr{H}^Q$ if and only if the
corresponding conditionally negative definite function
$g\mapsto\|v-g\cdot v\|^2$ factors through $Q$.\end{itemize}
\end{prop}
\begin{proof} 1) The first statement is immediate since, for all
$\lambda\in\R$, the function $(c,c')\mapsto(1-\lambda)c+\lambda
c'$ is continuous (actually 1-Lipschitz) on its domain of
definition.

2) If $\mathscr{H}$ is the Hilbert space of a unitary
representation, then $\mathscr{H}^Q$ is immediately seen to be a
linear subspace, and is closed by Proposition
\ref{prop:H^Q_closed}. Note that this also can be derived as a
particular case of 1). The nontrivial part of the last statement
in 2) follows from Proposition \ref{prop:H^Q_fct}.

3) is similar.\end{proof}

\begin{lem}
Let $G\to Q$ be a morphism with dense image, and $(\pi_i)$ is a
family of unitary representations of $G$. Then
$\left(\bigoplus\pi_i\right)^Q=\bigoplus\pi_i^Q$.\label{soussomme}
\end{lem}
\begin{proof} The inclusion $\bigoplus\pi_i^Q\subset
\left(\bigoplus\pi_i\right)^Q$ is trivial. Let $p_i$ denote the
natural projections, set $\pi=\bigoplus\pi_i$, and write
$\mu_\xi(g)=\pi(g)\xi$. Then, if $\xi$ is a $Q$-vector, i.e. if
$\mu_\xi$ factors through $Q$, then $p_i\circ\mu_\xi$ also factors
through $Q$. But $p_i\circ\mu_\xi=\mu_{p_i(\xi)}$, so that
$p_i(\xi)$ is a $Q$-vector.\end{proof}

\subsection{Resolutions}

\begin{conv}
If $G\to Q$ is a morphism with dense image, and $\pi$ is a unitary
representation of $G$ factoring through a representation
$\tilde{\pi}$ of $Q$, we write $1_Q\prec\pi$ rather than
$1_Q\prec\tilde{\pi}$ or $1\prec\tilde{\pi}$ to say that
$\tilde{\pi}$ almost has invariant vectors (note that $1_Q\prec
\pi$ implies $1_G\prec\pi$, but the converse is not true in
general). Similarly, if $(\pi_i)$ is a net of unitary
representations of $G$ factoring through representations
$\tilde{\pi}_i$ of $Q$, when we write $\pi_i\to 1_Q$, we mean for
the Fell topology on unitary representations of~$Q$.
\end{conv}

\begin{defn}[Resolutions]
Let $G$ be a locally compact group, and $f:G\to Q$ a morphism to
another locally compact group $Q$, such that $f(G)$ is dense
in~$Q$.

We say that $f$ is a {\em resolution} if, for every unitary
representation $\pi$ of $G$ almost having invariant vectors, then
$1_Q\prec \pi^Q$, meaning that $\pi^Q$, viewed as a representation
of $Q$, almost has invariant vectors (in particular, $\pi^Q\neq
0$).

We call $f$ a Haagerup resolution if $Q$ is Haagerup.
\end{defn}

The definition of resolution generalizes the notion of relative
Property~(T) of a closed normal subgroup $N$ of $G$, since $G\to
G/N$ is a resolution if and only if $(G,N)$ has relative
Property~(T).

\begin{rem}
When $G$ is $\sigma$-compact, it can be shown \cite[\S
2.3.8]{CorThe} that a morphism $f:G\to Q$ is a resolution if and
only if it satisfies the (a priori weaker) condition: every
unitary representation $\pi$ of $G$ such that $1_G\prec\pi$
satisfies $\pi^Q\neq 0$. The proof is not immediate and makes use
of the results in the sequel (namely, Theorem \ref{thm:resol}).
\end{rem}

\begin{prop}
Let $f:G\to Q$ be a morphism between locally compact groups, with
dense image. The following are equivalent:

\begin{itemize}\item[(1)] $f$ is a resolution.

\item[(2)] For every net $(\pi_i)$ of unitary representations  of
$G$ converging to $1_G$, we have $\pi_i^Q\to 1_Q$.\end{itemize}
\label{prop:resol_nets_1}\end{prop}

\begin{proof} (2)$\Rightarrow$(1) is trivial. Suppose (1). Let $\pi_i\to
1_G$. Then, for every subnet $(\pi_j)$,
$1_G\prec\bigoplus_j\pi_j$. By (1), $1_Q\prec
(\bigoplus_j\pi_j)^Q$, which equals $\bigoplus_j\pi_j^Q$ by Lemma
\ref{soussomme}. Hence $\pi_i^Q\to 1_Q$.\end{proof}

\begin{cor}
Let $G\to Q$ be a resolution. Then for every net $(\pi_i)$ of
irreducible unitary representations of $G$ converging to $1_G$,
eventually $\pi_i$ factors through a representation
$\tilde{\pi}_i$ of $Q$, and $\tilde{\pi}_i\to
1_Q$.\qed\label{cor:resol_irr}\end{cor}

The converse of Corollary \ref{cor:resol_irr} is more involved,
and is proved (Theorem \ref{thm:resol_irr}) under the mild
hypothesis that $G$ is $\sigma$-compact.

\medskip

Thus, a resolution allows to convey properties about the
neighbourhood of $1_Q$ in $\hat{Q}$ into properties about the
neighbourhood of $1_G$ in $\hat{G}$. For instance, this is
illustrated by Property $(\tau)$ (see Section
\ref{subs:application_tau}).

\medskip

The following proposition generalizes the fact that relative
Property~(T) is inherited by extensions (Proposition
\ref{prop:relT_extensions}), and is one of our main motivations
for having introduced resolutions.

\begin{prop}
Let $f:G\to Q$ be a resolution, and $X\subset G$. Then $(G,X)$ has
relative Property~(T) if and only if $(Q,f(X))$
does.\label{prop:resol_Trel}
\end{prop}
\begin{proof} The condition is trivially sufficient. Suppose that
$(Q,f(X))$ has relative Property~(T). Fix $\eps>0$, and let $\pi$
be a unitary representation of $G$ such that $1_G\prec\pi$. Since
$G\to Q$ is a resolution, $1_Q\prec\pi^Q$. Hence, by Property~(T),
$\pi^Q$ has a $(f(X),\eps)$-invariant vector; this is a
$(X,\eps)$-invariant vector for $\pi$.\end{proof}

Recall that a morphism between locally compact spaces is proper if
the inverse image of any compact subset is compact. It is easy to
check that a morphism $G\to H$ between locally compact groups is
proper if and only if its kernel $K$ is compact, its image $Q$ is
closed in $H$, and the induced map $G/K\to Q$ is an isomorphism of
topological groups.

\begin{cor}
Let $f:G\to Q$ be a Haagerup resolution. Then for every $X\subset
G$, $(G,X)$ has relative Property~(T) if and only if
$\overline{f(X)}$ is compact.

Accordingly, either $f$ is a proper morphism, so that $G$ is also
Haagerup, or there exists a noncompact closed subset $X\subset G$
such that $(G,X)$ has relative Property~(T). In particular, $G$
satisfies the TH alternative.\qed\label{cor:resol_Haag_TH_alt}
\end{cor}

\begin{thm}
Let $f:G\to Q$ be a resolution. Then $G$ is compactly generated if
and only if $Q$ is.\label{thm:resol_cpt_gen}
\end{thm}

This theorem generalizes compact generation of locally compact
groups with Property~(T)~\cite{K} (case when $Q=\{1\}$); in this
more specific direction, it generalizes Proposition 2.8 of
\cite{LZ}.

\begin{lem}
Let $f:G\to Q$ be a morphism with dense image between topological
groups, and let $\Omega$ be an open neighbourhood of 1 of $Q$.
Then, for all $n$, $f^{-1}(\Omega^n)\subset
f^{-1}(\Omega)^{n+1}$.\label{lem:mor_dense_preimages_iterees}
\end{lem}
\begin{proof} Let $x$ belong to $f^{-1}(\Omega^n)$. Write $f(x)=u_1\dots
u_n$ with $u_i\in\Omega$. Consider
$(\eps_1,\dots,\eps_n)\in\Omega^n$. Set $v_i=u_i\eps_i$ and
$v_0=u_1\dots u_n(v_1\dots v_n)^{-1}$, so that $f(x)=v_0v_1\dots
v_n$. If all $\eps_i$ are chosen sufficiently close to 1, then
$v_0\in\Omega$ and $u_i\eps_i\in\Omega$ for all $i$; by density of
$f(G)$, we can also impose that $u_i\eps_i\in f(G)$ for all~$i$.
We fix $\eps_1,\dots,\eps_n$ so that all these conditions are
satisfied. Since $v_0=f(x)(v_1\dots v_n)^{-1}$, we observe that
$v_0$ also belongs to $f(G)$. For all $i$, write $v_i=f(x_i)$, so
that $x=kx_0x_1\dots x_n$ with $k\in\textnormal{Ker}(f)$. Set
$y_0=kx_0$. Then $x=y_0x_1\dots x_n\in
f^{-1}(\Omega)^{n+1}$.\end{proof}

\begin{proof}[Proof of Theorem \ref{thm:resol_cpt_gen}] If $G$ is compactly
generated, so is $Q$ (only supposing that the morphism has dense
image). Indeed, if $K$ is a compact generating set of $G$, then
$f(K)$ generates a dense subgroup of $Q$. It follows that, if $K'$
is a compact subset of $Q$ containing $f(K)$ in its interior, then
$K'$ generates $Q$.

Conversely, suppose that $Q$ is compactly generated, and let
$\Omega$ be an open, relatively compact generating set. For every
open, compactly generated subgroup $H$ of $G$, let $\varphi_H$ be
its characteristic function. Then, when $H$ becomes big,
$\varphi_H$ converges to 1, uniformly on compact subsets of $G$.
Since $(G,f^{-1}(\Omega))$ has relative Property~(T) by
Proposition \ref{prop:resol_Trel}, it follows from Lemma
\ref{lem:Trel_inc_cpt_eng} that $f^{-1}(\Omega)$ is contained in a
compactly generated subgroup $H$ of $G$. By Lemma
\ref{lem:mor_dense_preimages_iterees}, $f^{-1}(\Omega^n)\subset
f^{-1}(\Omega)^{n+1}\subset H$. Since $Q=\bigcup\Omega^n$, it
follows that $G=H$.\end{proof}

\subsection{Lattices and resolutions}

The following theorem generalizes the fact that Property~(T) is
inherited by lattices.

\begin{thm}
Let $G$ be a locally compact group, $N$ a closed, normal subgroup.
Suppose that $(G,N)$ has relative Property~(T) (equivalently, the
projection $f:G\to G/N$ is a resolution).

Let $H$ be a closed subgroup of finite covolume in $G$, and write
$Q=\overline{f(H)}$. Then $f:H\to Q$ is a
resolution.\label{thm:MBL}
\end{thm}

Theorem \ref{thm:MBL} is a slight strengthening of \cite[Chap.
III, (6.3) Theorem]{Margulis}. To prove it, we need the following
lemma, whose ingredients are borrowed from \cite{BL}.

\begin{lem}
Let $G,N,H,Q$ be as in Theorem \ref{thm:MBL}. For every unitary
representation $\pi$ of $H$ factoring through $Q$, if
$1_H\prec\pi$, then $1_{Q}\prec\pi$. \label{BL_lem_almost-inv}
\end{lem}
\begin{proof} By \cite[Corollary 4.1(2)]{Jol}, $(f^{-1}(Q),N)$ has relative
Property~(T); hence, replacing $G$ by $f^{-1}(Q)$ if necessary, we
can suppose that $Q=G/N$.

Denote by $\tilde{\pi}$ the factorization of $\pi$ through $Q$,
and by $\hat{\pi}$ the (intermediate) factorization of $\pi$
through $G$, so that $\pi=\hat{\pi}|_H$.

By continuity of induction and since $H$ has finite covolume,
$1_G\prec\text{Ind}_H^G\pi$. On the other hand,
$\text{Ind}_H^G\pi=\text{Ind}_H^G\hat{\pi}|_H=\hat{\pi}\otimes
L^2(G/H)=\hat{\pi}\,\oplus\,(\hat{\pi}\otimes L^2_0(G/H))$.

We claim that $1_G\nprec\hat{\pi}\otimes L^2_0(G/H)$. It follows
that $1_G\prec\hat{\pi}$. Since every compact subset of $G/N$ is
the image of a compact subset of $G$ \cite[Lemma B.1.1]{BHV}, it
follows that $1_{G/N}\prec\tilde{\pi}$.

It remains to prove the claim. Since $\hat{\pi}|_N$ is a trivial
representation, $\hat{\pi}|_N\otimes L^2_0(G/H)|_N$ is a multiple
of $L^2_0(G/H)|_N$. But, by \cite[Lemma 2]{BL}, $L^2_0(G/H)$ does
not contain any nonzero $N$-invariant vector. Accordingly, neither
does $\hat{\pi}\otimes L^2_0(G/H)$. Hence, by relative
Property~(T), $1_G\nprec\hat{\pi}\otimes L^2_0(G/H)$.\end{proof}

\begin{proof}[Proof of Theorem \ref{thm:MBL}] Let $\pi$ be a unitary
representation of $H$, and suppose that $1_H\prec \pi$. Using
\cite[Chap. III, (6.3) Theorem]{Margulis} twice\footnote{The
assumption in \cite{Margulis} is that $N$ has Property~(T), but it
is clear from the proof that relative Property~(T) for $(G,N)$ is
sufficient.}, $\pi^Q\neq 0$, and its orthogonal in $\pi$ does not
almost contain invariant vectors. It follows that $1_H\prec\pi^Q$.
By Lemma \ref{BL_lem_almost-inv}, $1_Q\prec\pi^Q$.\end{proof}

\begin{rem}
The conclusion of Lemma \ref{BL_lem_almost-inv} is false if we
drop the assumption that $(G,N)$ has relative Property~(T), as the
following example shows.

Set $G=\Z\times\R/\Z$ and $N=\Z\times\{0\}$. Let $H$ be the cyclic
subgroup of $G$ generated by $(1,\alpha)$, where
$\alpha\in(\R-\mathbf{Q})/\Z$.

The projection $p:H\to\R/\Z$ has dense image. Hence, the
Pontryagin dual morphism: $\hat{p}:\Z\simeq \widehat{\R/\Z}\to
H^*\simeq\R/\Z$ also has dense image. Take a sequence $(\chi_n)$
of pairwise distinct nontrivial characters of $\R/\Z$ such that
$\hat{p}(\chi_n)$ tends to 0. Then the direct sum
$\pi=\bigoplus\chi_n$ does not weakly contain the trivial
representation (otherwise, since $\R/\Z$ has Property~(T), it
would contain the trivial representation), but $\pi\circ p|_H$
weakly contains the trivial representation
$1_H$.\label{rem:almost_inv_cercle}
\end{rem}

We can now combine the results of Section \ref{Sec_relT_lie_alg}
with Theorem \ref{thm:MBL}. Let $G$ be a finite direct product of
connected Lie groups and algebraic groups over local fields of
characteristic zero: $G=L\times\prod_{i=1}^nH_i(\K_i)$. Write
$R_T(G)=R_T(L)\times\prod_{i=1}^nR_T(H_i)(\K_i)$, where $R_T$ is
defined in Sections \ref{subs:relT_alg} and \ref{subs:relT_lie}.
Observe that, by Theorems \ref{thm:Trel_Lie} and \ref{thm:RT_alg},
$(G,R_T)$ has relative Property~(T) and $G/R_T$ is Haagerup.
Denote by $f:G\to G/R_T(G)$ the quotient morphism.

\begin{cor}
Let $G$ be a finite product of connected Lie groups and (rational
points of) algebraic groups over local fields of characteristic
zero. Let $\Gamma$ be a closed subgroup of finite covolume in $G$.
Then there exists a Haagerup resolution for $\Gamma$, given by
$f:\Gamma\to
\overline{f(\Gamma)}$.\qed\label{cor:resol_lattices_lie_alg}
\end{cor}

\subsection{The factorization Theorem}

\begin{thm}
Let $G$ be a locally compact group, $f:G\to Q$ a resolution,
$u:G\to H$ a morphism to a locally compact group, with dense
image, where $H$ is Haagerup. Then there exists a compact, normal
subgroup $K$ of $H$, and a factorization $Q\to H/K$ making the
following diagram commutative
$$\xymatrix{ G \ar[d]_u \ar[r]^f
& Q \ar@{.>}[d] \\
H \ar[r]^p & H/K. }$$\label{thm:factorization}
\end{thm}

\begin{proof} $H$ has a $C_0$ unitary representation $\pi$ with almost
invariant vectors. Therefore $\pi\circ u$ almost has invariant
vectors, so that, passing to a subrepresentation if necessary, we
can suppose that $\pi\circ u$ factors through a representation
$\tilde{\pi}$ of $Q$. We fix a normalized coefficient $\varphi$ of
$\tilde{\pi}\circ f$.

Let $(g_i)$ be a net in $G$ such that $f(g_i)\to 1$. Then
$\varphi(g_i)\to 1$. This implies that $(u(g_i))$ is bounded in
$H$, since $\varphi$ is a $C_0$ function. Let $K\subset H$ be the
set of all limits of $u(g_i)$ for such nets $(g_i)$. Then $K$ is a
compact, normal subgroup of $H$ (it is normal thanks to the
density of $u(G)$ in $H$).

Let $p$ be the projection: $H\to H/K$. We claim that $p\circ u$
factors through $Q$. Indeed, if $(g_i)$ is a net in $G$ such that
$(f(g_i))$ is Cauchy in $Q$, then $p\circ u$ is also Cauchy in
$H/K$. This implies that $p\circ u$ factors through
$Q$.\end{proof}

\subsection{Applications to Property ($\tau$) and related
properties}\label{subs:application_tau}

Recall that a representation of a group is said to be {\it finite}
if its kernel has finite index.

\begin{defn}
We recall that a topological group $G$ has Property ($\tau$)
[resp. ($\tau_{FD}$)] if for every net $(\pi_i)$ of finite (resp.
finite dimensional) irreducible unitary representations of $G$
converging to $1_G$, eventually $\pi_i=1_G$.

We say that a topological group $G$ has Property
($\textnormal{FH}_{FD}$) if every isometric action of $G$ on a
finite-dimensional Hilbert space has a fixed point. Equivalently,
every finite-dimensional unitary representation has vanishing
1-cohomology.

We say that a topological group $G$ has Property
($\textnormal{FH}_{F}$) if every finite unitary representation has
vanishing 1-cohomology.

The topological group $G$ has Property ($\textnormal{FAb}^\R$)
[resp. (FAb)] if for every closed subgroup of finite index $H$ of
$G$, we have $\textnormal{Hom}(H,\R)=0$ [resp.
$\textnormal{Hom}(H,\Z)=0$].
\end{defn}

It turns out that Properties ($\textnormal{FH}_{F}$) and
($\textnormal{FAb}^\R$) are equivalent. This is shown in
\cite{LZu}\footnote{The group there is assumed to be finitely
generated but this has no importance.} using induction of unitary
representations and Shapiro's Lemma. Alternatively, this can be
shown using induction of affine representations.

Note also that ($\textnormal{FAb}^\R$) implies (FAb), and they are
clearly equivalent for finitely generated groups; while $\R$ or
$\mathbf{Q}$ satisfy (FAb) but not ($\textnormal{FAb}^\R$).

Property $(\tau)$ clearly implies (FAb) \cite{LZu}; it is there
observed that the first Grigorchuk group does not have Property
$(\tau)$, but has Property ($\textnormal{FH}_{FD}$) since all its
linear representations are finite. However, no finitely presented
group is known to satisfy (FAb) but not Property $(\tau)$.

We provide below an example of a finitely generated group having
Property $(\tau)$ (hence ($\textnormal{FH}_{F}$)) but not
($\textnormal{FH}_{FD}$). We do not know if this is the first
known example.

We begin by a general result.

\begin{prop}
Let $G,Q$ be locally compact, and $G\to Q$ be a resolution. Let
(P) be one of the Properties: (T), ($\tau$), ($\tau_{FD}$), (FAb),
($\textnormal{FAb}^\R$), ($\textnormal{FH}_{F}$),
($\textnormal{FH}_{FD}$). %Let $G\to Q$ be a resolution, with $G,Q$ locally compact, $\sigma$-compact.
Then $G$ has Property (P) if and only if $Q$
does.\label{prop:resol_FAb_etc}
\end{prop}
\begin{proof} In all cases, Property (P) for $G$ clearly implies Property
(P) for $Q$.

Let us show the converse. For (T), ($\tau$), and ($\tau_{FD}$)
this follows directly from Proposition \ref{prop:resol_nets_1}.

Suppose that $G$ does not have Property (FAb). Let $N\subset G$ be
a closed normal subgroup of finite index such that
$\textnormal{Hom}(N,\Z)\neq 0$. Let $M$ be the kernel of a
morphism of $N$ onto $\Z$, and set $K=\bigcap_{g\in G/N}gMg^{-1}$.
Then $K$ is the kernel of the natural diagonal morphism
$N\to\prod_{g\in G/N}N/gMg^{-1}\simeq\Z^{G/N}$. It follows that
$N/K$ is a nontrivial free abelian group of finite rank, and $K$
is normal in $G$, so that $H=G/K$ is infinite, finitely generated,
virtually abelian. Since $H$ is Haagerup, by Theorem
\ref{thm:factorization}, $Q$ maps onto the quotient of $H$ by a
finite subgroup $F$. Since $H/F$ is also infinite, finitely
generated, virtually abelian, $Q$ does not have Property (FAb).

The case of Property ($\textnormal{FAb}^\R$) can be proved
similarly; since ($\textnormal{FAb}^\R$) is equivalent to
($\textnormal{FH}_{F}$) which is treated below, we omit the
details.

Suppose that $G$ does not have Property ($\textnormal{FH}_{FD}$).
Let $G$ act isometrically on a Euclidean space $E$ with unbounded
orbits, defining a morphism $\alpha:G\to\textnormal{Isom}(E)$. Set
$H=\overline{\alpha(G)}$. Since $\textnormal{Isom}(E)$ is
Haagerup, so is $H$. By Theorem \ref{thm:factorization}, $H$ has a
compact normal subgroup $K$ such that $Q$ has a morphism with
dense image into $H/K$. Observe that the set of $K$-fixed points
provides an action of $H/K$ on a non-empty affine subspace of $E$,
with unbounded orbits. So $Q$ does not have Property
($\textnormal{FH}_{FD}$).

The case of Property ($\textnormal{FH}_{F}$) can be treated
similarly, noting that $\alpha$ maps a subgroup of finite index to
translations, and this is preserved after restricting to the
action on an affine subspace.\end{proof}

\begin{rem}
The case of Property ($\textnormal{FH}_{FD}$) in Proposition
\ref{prop:resol_FAb_etc} contains as a particular case Theorem~B
in \cite{BL}, without making use of the Vershik-Karpushev Theorem,
while the proof given in \cite{BL} does.
\end{rem}

Proposition \ref{prop:resol_FAb_etc} justifies why we do not have
restricted the definitions of Property ($\tau$), etc., to discrete
groups (as is usually done), since in many cases, when we have a
resolution $G\to Q$, the group $Q$ is non-discrete. For instance,
all these properties are easy or trivial to characterize for
connected Lie groups.

\begin{prop}
Let $G$ be a connected locally compact group. Then

\begin{itemize}\item[1)] $G$ has Property $(\tau)$.

\item[2)] $G$ has Property ($\textnormal{FH}_{F}$) if and only if
$\textnormal{Hom}(G,\R)=0$.

\item[3)] $G$ has Property ($\tau_{FD}$) if and only if
$\textnormal{Hom}(G,\R)=0$.

\item[4)] $G$ has Property ($\textnormal{FH}_{FD}$) if and only if
every amenable quotient of $G$ is
compact.\end{itemize}\label{prop:tau_etc_lie}
\end{prop}
\begin{proof} Since $G$ has no proper closed finite index subgroup, 1) and
2) are immediate.

3) The condition is clearly necessary. Conversely, suppose that
$\textnormal{Hom}(G,\R)=0$. Let $W$ be the intersection of all
kernels of finite-dimensional unitary representations of $G$.
Clearly, it suffices to show that $G/W$ has Property
$(\tau_{FD})$. By \cite[Th\'eor\`eme 16.4.6]{Dix},
$G/W\simeq\R^n\times K$ for some compact group $K$. The assumption
then implies $n=0$, so that $G/W$ is compact, so that $G$ has
Property $(\tau_{FD})$.

4) Suppose that $G$ does not have Property
($\textnormal{FH}_{FD}$). Then there exists an unbounded isometric
affine action of $G$ on some Euclidean space. Let $K$ be a compact
normal subgroup of $G$ such that $G/K$ is a Lie group \cite{MZ}.
Restricting to the (non-empty) affine subspace of $K$-fixed
points, we can suppose that $K$ is contained in the kernel $N$ of
this action. Necessarily, the Lie group $G/N$ is not compact, and
amenable since it embeds in the amenable Lie group
$\textnormal{O}(n)\ltimes\R^n$.

Conversely, suppose $G$ has a noncompact amenable quotient $H$.
Since $H$ does not have Property~(T), by a result of Shalom
\cite{Sha00} (see \cite[Section 3.2]{BHV}), there exists an
irreducible unitary representation $\pi$ of $H$ with non-vanishing
1-reduced cohomology. By \cite[Theorem~3.1]{MartinJLT}, $\pi$ is
finite-dimensional\footnote{It is possible to prove 4) more
directly, but we have used Shalom's and Martin's results to make
short.}.\end{proof}

\begin{prop}
Fix $n\ge 5$, set $\Gamma=\SO_n(\Z[2^{1/3}])\ltimes
\Z[2^{1/3}]^n$. Then $\Gamma$ is a finitely presentable group, has
Property $(\tau_{FD})$ (hence Property $(\tau)$, hence Property
($\textnormal{FH}_{F}$)), but not
($\textnormal{FH}_{FD}$).\label{prop:SO5_ltimes_R5}
\end{prop}
\begin{proof} Note that $\SO_n(\C)$ and $\SO_n(\C)\ltimes\C^n$ have
Property~(T). Since $\Gamma$ is an irreducible lattice in the
connected Lie group
$(\SO_n(\R)\ltimes\R^n)\times(\SO_n(\C)\ltimes\C^n)$, it is
finitely presentable, and, moreover, by Theorem \ref{thm:MBL},
$\Gamma\to\SO_n(\R)\ltimes\R^n$ is a resolution.

By Proposition \ref{prop:tau_etc_lie}, $\SO_n(\R)\ltimes\R^n$ has
Property $(\tau_{FD})$. Thus $\Gamma$ also has Property
$(\tau_{FD})$ by Proposition \ref{prop:resol_FAb_etc}. On the
other hand, the embedding of $\Gamma$ in $\SO_n(\R)\ltimes\R^n$
provides an isometric action of $\Gamma$ with unbounded orbits on
the $n$-dimensional Euclidean space.\end{proof}

\begin{rem}
It is asked in \cite{LZ} whether there exists a finitely generated
group with Property $(\tau)$ but not $(\tau_{FD})$. Obvious
non-finitely generated examples are $\mathbf{Q}$ and $\R$. It may
be tempting to find a finitely generated group $\Gamma$ with a
resolution $\Gamma\to\R$, but unfortunately no such $\Gamma$
exists. Indeed, since $\Gamma$ is discrete and
$\textnormal{Hom}(\Gamma,\R)\neq 0$, there exists a discrete,
nontrivial, torsion-free abelian quotient $\Lambda$ of $\Gamma$.
By Theorem \ref{thm:factorization}, there exists a factorization:
$\R\to\Lambda$, necessarily surjective. This is a contradiction
since $\R$ is connected. Thus the question remains open.
\end{rem}

%%%%%%%%%%%%%%%%%%%%%%%%%%%%%%%%%%%%%%%%%%%%%%%%%%%%%%%%%%%%%%%%%%%%%%%%%%%%%%%%%

\subsection{Subgroups of simple Lie groups}

Let $G$ be a connected simple Lie group, with Lie algebra $\g$. We
are interested in subgroups $\Gamma\subset G$, viewed as discrete
groups.

1) If $\g$ is isomorphic to
$\mk{sl}_2(\R)\simeq\mk{so}(2,1)\simeq\mk{su}(1,1)$,
$\mk{sl}_2(\C)\simeq\mk{so}(3,1)$, or $\mk{so}_3(\R)$, then every
$\Gamma\subset G$ is Haagerup \cite[\S 5, Theorem~4]{GHW} (see
\cite{CorJLT} for the case of $\widetilde{\SL_2(\R)}$).

2) If $G$ has Property~(T) and
$\g\simeq\!\!\!\!\!/\;\mk{so}_3(\R)$, then there exists
$\Gamma\subset G$ with Property~(T): if $G$ is noncompact, take
any lattice. If $G$ is compact, this is due to Margulis
\cite[chap. III, Proposition 5.7]{Margulis}.

3) If $\g\simeq\mk{so}(n,1)$ with $n\ge 5$ or
$\g\simeq\mk{su}(n,1)$ with $n\ge 3$, then then there exists an
infinite subgroup $\Gamma\subset G$ with Property~(T): it suffices
to observe that $G$ contains a subgroup locally isomorphic to
$\SO(n)$ ($n\ge 5$) or $\SU(n)$ ($n\ge 3$), and such a subgroup
contains an infinite subgroup with Property~(T) by 2).

4) There are only two remaining cases: $\g\simeq\mk{so}(4,1)$ and
$\g\simeq\mk{su}(2,1)$. We are going to show that the behaviour
there is different from that in preceding examples.

The only result already known is that if $\g\simeq\mk{so}(4,1)$ or
$\g\simeq\mk{su}(2,1)$, then no infinite $\Gamma\subset G$ can
have Property~(T); this follows from 1) since such $\Gamma$ would
be contained in a maximal compact subgroup. This result is
generalized in the following theorem.

\begin{thm}
Let $G$ be a connected Lie group, locally isomorphic to either
$\SO(4,1)$ or $\SU(2,1)$. Let $\Gamma\subset G$ be any subgroup,
and view $\Gamma$ as a discrete group.

\begin{itemize}\item[1)] If $\Lambda\subset\Gamma$ is a normal subgroup such that
$(\Gamma,\Lambda)$ has relative Property~(T), then $\Lambda$ is a
finite subgroup of $G$.

\item[2)] If $\Gamma$ is not dense, and if $\Lambda\subset\Gamma$
is a subgroup such that $(\Gamma,\Lambda)$ has relative
Property~(T), then $\Lambda$ is a finite subgroup of $G$.

\item[3)] If $\Gamma$ is dense, and $X\subset\Gamma$ is a normal
subset (i.e. invariant under conjugation) such that $(\Gamma,X)$
has relative Property~(T), then $X$ is a finite subset of the
centre of $G$.\label{thm:sgSO41SU21}\end{itemize}

Suppose that $G$ is locally isomorphic to $\SU(2,1)$. Then we have
stronger statements:

\begin{itemize}\item[4)] If $\Gamma$ is not dense, then $\Gamma$ is Haagerup.

\item[5)] If $X\subset\Gamma$ is a normal subset and $(\Gamma,X)$
has relative Property~(T), then $X$ is a finite subset of
$G$.\end{itemize}\label{thm:relT_SO41_SU21}\end{thm}
\begin{proof} Fix a
subset $X\subset\Gamma$ such that $(\Gamma,X)$ has relative
Property~(T). We make a series of observations.

a) First note that $G$ is Haagerup, a fact due to \cite{FH} (and
\cite[Chap. 4]{CCJJV} in the case of $\widetilde{\SU}(2,1)$).
Therefore, by relative Property~(T), $\overline{X}$ must compact.
Denote by $\mk{h}$ the Lie algebra of $\overline{\Gamma}$.

b) Suppose, in this paragraph b), that $\Gamma$ is dense in $G$,
i.e. $\mk{h}=\g$; and suppose that $X$ is a normal subset. Then
$\overline{X}$ is a compact, normal subset in $G$. Let $Z$ be the
centre of $G$, and fix $h\in X$. Then the conjugacy class of $h$
in $G/Z$ is relatively compact. Let $M$ be the symmetric space
associated to $G/Z$, and fix $y\in M$. Then the function $g\mapsto
d(ghg^{-1}y,y)$ is bounded, so that $h$ has bounded displacement
length. Since $M$ is CAT($-1$) and geodesically complete, this
implies that $h$ acts as the identity, i.e. $h\in Z$. Accordingly,
$X\subset Z$, so that $X$ is discrete. Since it is relatively
compact, it is finite. This proves 3).

c) Suppose that $\Gamma$ is Zariski dense modulo $Z$, but not
dense. Then $\Gamma$ is discrete. Indeed, $\Gamma$ is contained in
the stabilizer $W$ of $\mk{h}$ for the adjoint action. Since $W$
is Zariski closed modulo $Z$, this implies that $\mk{h}$ is an
ideal in $\g$, so that, since $\mk{h}\neq\g$ and $\g$ is simple,
$\mk{h}=\{0\}$, i.e. $\Gamma$ is discrete. Since $G$ is Haagerup,
this implies that $\Gamma$ is Haagerup.

d) Now suppose that $Z=1$ and $\Gamma$ is not Zariski dense. Let
$N$ be the Zariski closure of $\Gamma$. Let $R_u$ be the unipotent
radical of $N$, and $L=CS$ a Levi factor, with abelian part $C$,
and semisimple part $S$. The possibilities for simple factors in
$S$ are rather restricted. The complexification $G_\C$ is
isomorphic to either $\textnormal{PSO}_5(\C)$ or
$\textnormal{PSL}_3(\C)$. In both cases, by a dimension argument,
the only possible simple subgroups of $G_\C$ are, up to isogeny,
$\SL_2(\C)$, and maybe $\SL_3(\C)$ in $\textnormal{PSO}_5(\C)$;
however, $\mk{sl}_3(\C)$ does not embed in $\mk{so}_5(\C)$ as we
see, for instance, by looking at their root systems. So the only
possible factors in $S$ are, up to isogeny, $\SL_2(\C)$,
$\SL_2(\R)$, and $\SO_3(\R)$. By \cite[\S 5, Theorem~4]{GHW}, the
image of $\Gamma$ in $N/R_u$ is Haagerup, so that the image of $X$
in $H/R_u$ is finite.

e) We keep the assumptions of d), and suppose moreover that
$X=\Lambda$ is a subgroup. Since $\Lambda$ is relatively compact,
and $R_u$ is unipotent, $\Lambda\cap R_u=\{1\}$. Since we proved
in d) that the image of $\Lambda$ in $N/R_u$ is finite, this
implies that $\Lambda$ is finite.

Now let us drop the assumption $Z=1$. Then the image of $\Lambda$
modulo $Z$ is finite, so that, by the case $Z=1$, $\Lambda$ is
virtually contained in $Z$. This implies that $\Lambda$ is
discrete, hence finite since it is also relatively compact.

In view of c), d), and e), 2) is now proved; observe that 1) is an
immediate consequence of 2) and 3).

f) Now suppose that $\g\simeq\mk{su}(2,1)$, and let us prove 4).
Observe that 5) is an immediate consequence of 3) and 4).

We first suppose that $Z=1$, and that $\Gamma$ is not Zariski
dense. So we continue with the notation of d). Write $S=S_cS_{nc}$
by separating compact and noncompact simple factors.

Suppose that $S_{c}\neq 1$. This is a compact subgroup, up to
conjugate, we can suppose that it is contained in the maximal
subgroup $\text{PS}(\text{U}(2)\times\text{U}(1))$. The Lie
algebra of $S_c$ is identified with $\mk{su}(2)$.

\begin{cla}The only proper subalgebra of $\mk{su}(2,1)$ properly containing
$\mk{su}(2)$ is $\mk{s}(\mk{u}(2)\times\mk{u}(1))$.
\end{cla}
Let us prove the claim. Let $\mk{k}$ be such a subalgebra. If
$\mk{k}\subset \mk{s}(\mk{u}(2)\times\mk{u}(1))$, then
$\mk{k}=\mk{s}(\mk{u}(2)\times\mk{u}(1))$ by a dimension argument.

Otherwise, we claim that the action of $\mk{k}$ on $\C^3$ is
irreducible. Let us consider the decomposition
$\C^3=\C^2\oplus\C$. Let $A$ be the $\C$-subalgebra of
$\text{M}_3(\C)$ generated by $\mk{k}$. Since the action of
$\mk{su}(2)$ on $\C^2$ is irreducible, $A$ contains
$\text{M}_2(\C)\times\text{M}_1(\C)$. In particular, the only
possible stable subspaces are $\C^2\oplus\{0\}$ and
$\{0\}\oplus\{0\}\oplus\C$. Now observe that since they are
orthogonal to each other, if one is stable by $\mk{k}$, then so is
the other. So, if $\mk{k}$ does not act irreducibly, it preserves
these subspace; this means that
$\mk{k}\subset\mk{s}(\mk{u}(2)\times\mk{u}(1))$. This proves the
claim.

\medskip

By the claim, $N$ is virtually isomorphic to a connected Lie group
locally isomorphic to either $\SO_3(\R)$ or $\SO_3(\R)\times\R$.
So, by \cite[\S 5, Theorem~4]{GHW}, $\Gamma$ is Haagerup.

Otherwise, $S_c=1$. Since $G$ is Haagerup, by \cite[Chap.
4]{CCJJV} (or Corollary \ref{cor:Lie_haag}), $[S_{nc},R_u]=\{1\}$
so that $S_{nc}$ is, up to a finite kernel, a direct factor of
$N$. Since we have proved in d) that the only possible simple
factor appearing in $S_{nc}$ are locally isomorphic to $\SL_2(\R)$
or $\SL_2(\C)$,\footnote{Actually, it is easily checked that
$\mk{sl}_2(\C)$ does not embed in $\mk{su}(2,1)$.} in view of
\cite[\S 5, Theorem~4]{GHW}, this implies that $\Gamma$ has a
subgroup of finite index having the Haagerup Property, so that
$\Gamma$ is Haagerup.

Finally let us drop the hypothesis $Z=1$. Let $N$ be the preimage
in $G$ of the Zariski closure of $\Gamma$ in $G/Z$. There are two
possible cases:

$\bullet$ $N$ has finitely many connected components. Then, by
\cite[Theorem~3.13]{CorJLT} (which relies on similar arguments),
every subgroup of $N$ is Haagerup for the discrete topology.

$\bullet$ $N$ has infinitely many connected components. Then $N$
is almost the direct product of $Z$ and $N/Z$, so that, by the
case $Z=1$, every subgroup of $N$ is Haagerup for the discrete
topology.\end{proof}

The following proposition shows that the statements in Theorem
\ref{thm:relT_SO41_SU21} are, in a certain sense, optimal: in 1),
the assumption that $\Lambda$ be a normal subgroup cannot be
dropped, etc.

\begin{prop}
Let $G$ be a connected Lie group, locally isomorphic to either
$\SO(4,1)$ or $\SU(2,1)$.

\begin{itemize}\item[(1)] $G$ has finitely presented subgroups
$\Gamma\supset\Lambda$, such that $\Lambda$ is infinite and
$(\Gamma,\Lambda)$ has relative Property~(T).

\item[(2)] If $G$ is locally isomorphic to $\SO(4,1)$, then $G$
has a finitely presented subgroup $\Gamma$ and an infinite normal
subset $X\subset\Gamma$ such that $(\Gamma,X)$ has relative
Property~(T).\end{itemize}\label{prop:sg_trel_SO41_SU21}
\end{prop}

\begin{proof} (1) First suppose that $G=\SU(2,1)$, and write $G=H(\R)$,
where $H(R)$ is defined, for every commutative ring $R$ as the set
of matrices $(A,B)$ (rather denoted $A+iB$) satisfying the
relation $(\,^t\!A-i\,^t\!B)J(A+iB)=J$, where\footnote{This
relation must be understood as a relation where $i$ is a formal
variable satisfying $i^2=-1$. In other words, this means
$\,^t\!AJA+\,^t\!BJB=J$ and $\,^t\!AJB-\,^t\!BJA=0$.} $J$ is the
diagonal matrix $\textnormal{diag}(1,1,-1)$. Let $K$ be defined as
the upper-left $2\times 2$ block in $H$, so that
$K(\R)\simeq\SU(2)$. Observe that $H(\C)\simeq\SL_3(\C)$ and
$K(\C)\simeq\SL_2(\C)$.

Then $\Gamma=H(\Z[2^{1/3}])$ embeds as a lattice in $H(\R)\times
H(\C)$. By Theorem \ref{thm:MBL}, the projection $p$ of $\Gamma$
into $H(\R)\simeq\SU(2,1)$ is a resolution. Set
$\Lambda=K(\Z[2^{1/3}])$. Then $\Lambda$ is a lattice in
$K(\R)\times K(\C)$, so embeds as a cocompact lattice in
$K(\C)\simeq\SL_2(\C)$. On the other hand, since $p(\Lambda)$ is
relatively compact (it is dense in $\SU(2)$), by Proposition
\ref{prop:resol_Trel}, $(\Gamma,\Lambda)$ has relative
Property~(T). Note that, as lattices in connected Lie groups, they
are finitely presentable.

Let us now suppose that $G$ is locally isomorphic to $\SU(2,1)$,
and let $Z$ be its centre. Let $\Gamma,\Lambda$ be as above, and
let $\Gamma_0,\Lambda_0$ be their projection in
$G/Z\times\SL_3(\C)$. Finally, let $\Gamma_1,\Lambda_1$ be their
preimage in $G\times\SL_3(\C)$. If $Z$ is finite, then it is
immediate that $(\Gamma_1,\Lambda_1)$ has relative Property~(T),
and that they are finitely presented. So we suppose that
$G=\widetilde{\SU(2,1)}$. Let $q$ be the projection $G\times
\SL_3(\C)\to \SU(2,1)\times\SL_2(\C)$, and observe that
$\Gamma_1=q^{-1}(\Gamma)$ and $\Lambda_1=q^{-1}(\Lambda)$.

Since $K(\R)\times K(\C)\simeq\SU(2)\times\SL_2(\C)$ is simply
connected, $W=q^{-1}(K(\R)\times K(\C))$ is isomorphic to
$K(\R)\times K(\C)\times\Z$, and contains $\Lambda_1$ as a
lattice. So we can define $\Lambda_2$ as the projection of
$\Lambda_1$ into the unit component $W_0$, which is isomorphic to
$\Lambda$, hence finitely presentable. Since the projection of
$\Lambda_2$ on $G=K(\R)$ is relatively compact, by Proposition
\ref{prop:resol_Trel}, $(\Gamma_1,\Lambda_2)$ has relative
Property~(T).

\medskip

A similar example can be constructed in $\SO(4,1)$, projecting an
irreducible lattice from $\SO(4,1)\times\SO_5(\C)$. Since
$\SO(4,1)$ has finite fundamental group, we do not have to care
with some of the complications of the previous example.

(2) Observe that $\SO(4,1)$ has a subgroup isomorphic to
$\SO_3(\R)\ltimes\R^3$. Indeed, if we write $\SO(4,1)$ as $\{A|\;
\,^t\!AJA=J,\,\det(A)=1\}$, where $J=\begin{pmatrix}
  0 & 0 & 1 \\
  0 & I_3 & 0 \\
  1 & 0 & 0 \\
\end{pmatrix}$, then it contains the following subgroup, which is isomorphic to
$\SO_3(\R)\ltimes\R^3$:
$$P=\left\{\begin{pmatrix}
  1 & -\,^t\!vA & -\,^t\!vv/2 \\
  0 & A & v \\
  0 & 0 & 1 \\
\end{pmatrix}|\;A\in \SO(3),\; v\in \R^3\right\}.$$

Now consider the subgroup
$\Gamma=\SO_3(\Z[2^{1/3}])\ltimes\Z[2^{1/3}]^3$. Then $\Gamma$
embeds as a lattice in $(\SO_3(\C)\ltimes\C^3)\times
(\SO_3(\R)\ltimes\R^3)$. By Theorem \ref{thm:Trel_Lie},
$((\SO_3(\C)\ltimes\C^3)\times (\SO_3(\R)\ltimes\R^3),\C^3)$ has
relative Property~(T). Therefore, by Theorem \ref{thm:MBL}, the
inclusion morphism $\Gamma\to\SO_3(\Z[2^{1/3}])\ltimes\R^3$ is a
resolution. Let $B$ be the Euclidean unit ball in $\R^3$. Then, by
Proposition \ref{prop:resol_Trel}, $(\Gamma,\Z[2^{1/3}]^3\cap B)$
has relative Property~(T). Finally observe that
$X=\Z[2^{1/3}]^3\cap B$ is a normal subset in $\Gamma$.

Now observe that $\Gamma$ is contained in $P$, hence is contained
in the unit component $\SO_0(4,1)$. The only other connected Lie
group with Lie algebra $\mk{so}(4,1)$ is its universal covering
(of degree 2); taking the preimage of $\Gamma$ and $X$, we obtain
the required pair with relative Property~(T).\end{proof}

\begin{rem}
Examples similar to
$\Gamma=\SO_3(\Z[2^{1/3}])\ltimes\Z[2^{1/3}]^3$ (see the proof of
Proposition \ref{prop:sg_trel_SO41_SU21}) were already introduced
in \cite{CorJLT}. It was observed there that they provide the
first known examples of groups without the Haagerup Property
having no infinite subgroup with relative Property~(T). We have
made here more concrete the negation of the Haagerup Property by
exhibiting an infinite subset with relative
Property~(T).\label{rem:SO3_ltimes_R3}
\end{rem}

%%%%%%%%%%%%%%%%%%%%%%%%%%%%%%%%%%%%%%%%%%%%%%%%%%%%%%%%%%%%%%%%%%%%%%%%%%%%%%%%%%

\subsection{Affine resolutions}

Although they are probably known to the specialists, we found no
reference for the following lemmas.

\begin{lem}
Let $M$ be a complete CAT(0) metric space. Let $X$ be a nonempty
bounded subset, and let $B'(c,r)$ be the closed ball of minimal
radius containing $X$ \cite[Chap. II, Corollary 2.8(1)]{BH}.
Suppose that $X$ is contained in another ball $B'(c',r')$. Then
$$d(c,c')^2\le r'^2-r^2.$$\label{lem_centre}\end{lem}
\begin{proof} Set $d=d(c,c')$. Suppose the contrary, so that $d^2>r'^2-r^2\ge 0$.
For $t\in [0,1]$, set $p_t=(1-t)c+tc'$, which is a well-defined
point on the geodesic segment $[cc']$.

By \cite[(**) p.153]{Brown}, for every $z\in M$ and every $t\in
[0,1]$, \begin{equation}d(z,p_t)^2\le
(1-t)d(z,c)^2+td(z,c')^2-t(1-t)d^2.\label{eq:CAT0}\end{equation}
It follows that if $z\in B'(c,r)\cap B'(c',r')$, then
$d(z,p_t)^2\le (1-t)r^2+tr'^2-t(1-t)d^2$; denote by $u(t)$ this
expression. By an immediate calculation, $u(t)$ is minimal for
$t=t_0=(d^2+r^2-r'^2)/(2d^2)$, which belongs to $]0,1]$ by
assumption. Since $u(0)=r^2$, it follows that $u(t_0)<r^2$. Since
this is true for all $z\in B'(c,r)\cap B'(c',r')$, this implies
that $X$ is contained in a closed ball of radius $u(t_0)^{1/2}<r$,
contradiction.\end{proof}

\begin{lem}
Let $M$ be a complete CAT(0) metric space. Let $K$ be a nonempty
closed convex bounded subset, and $B'(c,r)$ the ball of minimal
radius containing $K$. Then $c\in
K$.\label{lem:convex_cont_centre}
\end{lem}
\begin{proof} Suppose that $c\notin K$, and let $p$ be its projection on
$K$ \cite[Chap. II, Proposition 2.4]{BH}. Fix $x\in K$. Then for
every $p'\in [px]$, $d(p,c)\le d(p',c)$. Hence, by
(\ref{eq:CAT0}), for all $t\in [0,1]$, $d(p,c)^2\le
(1-t)d(p,c)^2+td(x,c)^2-t(1-t)d(p,x)^2$. Taking the limit, after
dividing by $t$, when $t\to 0$, gives $d(x,p)^2\le
d(x,c)^2-d(p,c)^2$, so that $d(x,p)^2\le r^2-d(p,c)^2$. In other
words, $K\subset B'(p,(r^2-d(p,c)^2)^{1/2})$. This contradicts the
minimality of $r$.\end{proof}

\begin{lem}
Let $M$ be a complete CAT(0) metric space. Let $(F_n)$ be a
decreasing sequence of nonempty closed convex bounded subsets.
Then $\bigcap F_n\neq\emptyset$.\label{lem:intersection_convex}
\end{lem}
\begin{proof} Let $B'(c_n,r_n)$ be the ball of minimal radius containing
$F_n$. Observe that $c_n\in F_n$ by Lemma
\ref{lem:convex_cont_centre}. Moreover, $(r_n)$ is non-increasing,
hence converges.

On the other hand, if $m\le n$, then $F_n\subset F_m$. Applying
Lemma \ref{lem_centre}, we get $d(c_n,c_m)^2\le r_n^2-r_m^2$.
Therefore, $(c_n)$ is Cauchy, hence has a limit $c$, which belongs
to $\bigcap F_n$.\end{proof}

\begin{thm}
Let $f:G\to Q$ be a morphism with dense image between locally
compact groups. Let $G$ act by isometries on a complete CAT(0)
metric space $M$. Suppose that there exists a neighbourhood
$\Omega$ of $1$ in $Q$, such that, for some $w\in M$,
$f^{-1}(\Omega)w$ is bounded. Then $M^Q$ is nonempty.
\label{thm:X^Q_non_vide}
\end{thm}
\begin{proof} Let $(\Omega_n)$ be a sequence of compact symmetric
neighbourhoods of 1 in $Q$, contained in $\Omega$, such that
$\Omega_{n+1}\cdot\Omega_{n+1}\subset\Omega_n$ for all $n$. Set
$V_n=f^{-1}(\Omega_n)$.

By the assumption on $\Omega$, $V_n\cdot w$ is bounded for all
$n$. Let $B'(c_n,r_n)$ be the minimal ball containing $V_n\cdot
w$. Note that the sequence $(c_n)$ is bounded since $d(c_n,w)\le
r_n\le r_0$ for all $n$.

Then, for all $g\in V_{n+1}$, we have $g^{-1}V_{n+1}\cdot w\subset
V_n\cdot w\subset B'(c_n,r_n)$, so that $V_{n+1}\cdot w\subset
B'(gc_n,r_n)$. By Lemma \ref{lem_centre}, we have
$$d(c_{n+1},gc_n)\le \sqrt{(r_n-r_{n+1})(r_n+r_{n+1})}\le\sqrt{2r_0(r_n-r_{n+1})}.$$

Specializing this inequality to $g=1$, we obtain

$$d(c_{n+1},c_n)\le\sqrt{2r_0(r_n-r_{n+1})},$$
 and combining the two previous inequalities, we get, for all
 $g\in V_{n+1}$,

$$d(c_n,gc_n)\le d(c_{n+1},c_n)+d(c_{n+1},gc_n)\le 2\sqrt{2r_0(r_n-r_{n+1})}.$$

Set $u(n)=\sup\{2\sqrt{2r_0(r_m-r_{m+1})}|\;m\ge n\}$. Since
$(r_n)$ is non-increasing and nonnegative, $r_n-r_{n+1}\to 0$, so
that $u(n)\to 0$.

Note that, for all $g$, the function $x\mapsto d(x,gx)$ is
continuous and convex on $M$ \cite[Chap. II, Proposition 6.2]{BH}.
It follows that $F_n=\{v\in\HH|\;\forall g\in V_{n+1},\;d(v,gv)\le
u(n)\}$ is closed and convex. Set $K_n=F_n\cap B'(w,r_0)$. Then
$(K_n)$ is a decreasing sequence of closed, convex, bounded
subsets of $M$, nonempty since $c_n\in K_n$. By Lemma
\ref{lem:intersection_convex}, $\bigcap K_n$ is nonempty; pick a
point $y$ in the intersection. We claim that $y\in X^Q$: to see
this, let us appeal to Proposition \ref{prop:H^Q_fct}. Let $g_i$
be a net in $G$ such that $f(g_i)\to 1$. Set $n_i=\sup\{n|g_i\in
\Omega_{n+1}\}\in\N\cup\{\infty\}$. Then $n_i\to\infty$ since all
$\Omega_n$ are neighbourhoods of $1$ in $Q$, and $d(y,g_iy)\le
u(n_i)$ for all $i$ (where we set $u(\infty)=0$). It follows that
$d(y,g_iy)\to 0$. By Proposition \ref{prop:H^Q_fct}, $y\in
X^Q$.\end{proof}

\begin{defn}Let $f:G\to Q$ be a morphism with dense image between locally
compact groups. We call it an affine resolution if, for every
isometric action of $G$ on an affine Hilbert space, there exists a
non-empty $G$-invariant affine subspace such that the action of
$G$ on this subspace factors through $Q$.
\end{defn}

\begin{thm}

Let $G,Q$ be locally compact groups, $f:G\to Q$ be a morphism with
dense image. Consider the following conditions:

\begin{itemize}

\item[(1)] $G\to Q$ is a resolution.

\item[(2)] $(G,X)$ has relative Property~(T) for all subsets
$X\subset G$ such that $\overline{f(X)}$ is compact.

\item[(3)] $(G,X)$ has relative Property (FH) for all subsets
$X\subset G$ such that $\overline{f(X)}$ is compact.

\item[(4)] $G\to Q$ is an affine resolution.

\end{itemize}

Then the implications
(1)$\Rightarrow$(2)$\Rightarrow$(3)$\Leftrightarrow$(4) hold.
Moreover, if $G$ is $\sigma$-compact, then (4)$\Rightarrow$(1), so
that they are all equivalent.\label{thm:resol}\end{thm}
\begin{proof} (3) is
an immediate consequence of (4). The converse actually follows
immediately from Theorem \ref{thm:X^Q_non_vide}. The implication
(1)$\Rightarrow$(2) has been proved in Proposition
\ref{prop:resol_Trel}, and (2)$\Rightarrow$(3) follows from
Theorem \ref{thm:Del-Gui_relT}. It remains to prove
(4)$\Rightarrow$(1). Hence, suppose that $G$ (hence $Q$) is
$\sigma$-compact, and that $G\to Q$ is an affine resolution.

\begin{cla}
For every unitary representation $\pi$ of $G$ such that
$1_G\prec\pi$, $\pi^Q\neq 0$.
\end{cla}
\noindent Let us prove the claim. Let $\pi$ be a unitary
representation of $G$ on a Hilbert space $\HH$, such that
$1_G\prec\pi$. We must show that $\pi^Q\neq 0$. If $1\le\pi$, this
is trivially satisfied. So we can suppose that $1\nleqslant\pi$.
By a result of Guichardet which uses $\sigma$-compactness, (see
\cite[Theorem~2.13.2]{BHV}), $B^1(G,\pi)$ is not closed in
$Z^1(G,\pi)$, so that, in particular, $H^1(G,\pi)\neq 0$. Consider
$b\in Z^1(G,\pi)-B^1(G,\pi)$, and let $\alpha$ be the associated
affine action. Since $f$ is an affine resolution, $\alpha^Q$ is a
nonempty closed affine subspace $V$ of $\HH$. Then $V$ is not
reduced to a point $\{v\}$: otherwise, $v$ would be a fixed point
for the action of $G$, contradicting $b\notin B^1(G,\pi)$. Hence
the linear part of $\alpha^Q$ is a nonzero subrepresentation of
$\pi$, so that $\pi^Q$ is nonzero. This proves the claim.

\medskip

Let $\pi$ be a unitary representation of $G$ on a Hilbert space
$\HH$, such that $1_G\prec\pi$. We must show that $1_Q\prec\pi^Q$.
Again, since the case when $1\le\pi$ is trivial, we suppose that
$1\nleqslant\pi$. Let $\rho$ be the orthogonal of $\pi^Q$. By the
claim, $1_G\nprec\rho$. It follows that $1_G\prec\pi^Q$, so that
we can suppose that $\pi=\pi^Q$, i.e. $\pi$ factors through a
representation $\tilde{\pi}$ of $Q$.

\begin{cla}
The natural continuous morphism $\hat{f}:Z^1(Q,\tilde{\pi})\to
Z^1(G,\pi)$ is bijective.
\end{cla}
Let us prove the claim. The morphism $\hat{f}$ is clearly
injective. Take $b\in Z^1(G,\pi)$. Since $f$ is an affine
resolution, one can write $b(g)=b'(g)+\pi(g)v-v$ ($\forall g\in
G$), where $b'\in Z^1(G,\pi)$ factors through $Q$ and $v\in\HH$.
Since $\pi$ also factors through $Q$, this implies that $b$ does
so, which means that $b$ belongs to $\textnormal{Im}(\hat{f})$,
and the claim is proved.

\medskip

Since $G$ and $Q$ are $\sigma$-compact, $Z^1(G,\pi)$ and
$Z^1(Q,\tilde{\pi})$ are Fr\'echet spaces. Since
$\hat{f}:Z^1(Q,\tilde{\pi})\to Z^1(G,\pi)$ is bijective, by the
open mapping Theorem, it is an isomorphism. Note that it maps
$B^1(Q,\tilde{\pi})$ bijectively onto $B^1(G,\pi)$, and that
$B^1(G,\pi)$ is not closed in $Z^1(G,\pi)$, as we used in the
proof of the first claim. It follows that $B^1(Q,\tilde{\pi})$ is
not closed in $Z^1(Q,\tilde{\pi})$. Using the converse in
Guichardet's result \cite[Theorem 2.13.2]{BHV},
$1_Q\prec\tilde{\pi}$.\end{proof}

As a corollary of Theorems \ref{thm:X^Q_non_vide} and
\ref{thm:resol}, we get:

\begin{cor}
Let $G,Q$ be locally compact groups, and let $f:G\to Q$ be a
resolution. Let $G$ act isometrically on a complete CAT(0) metric
space $M$, and suppose that there exists a $G$-equivariant proper
embedding $i$ of $M$ in a Hilbert space. Then
$M^Q\neq\emptyset$.\label{cor:resol_met_hilb}
\end{cor}
\begin{proof} Let $\Omega$ be a compact neighbourhood of $1$ in $Q$, and
set $V=f^{-1}(\Omega)$. Fix $x\in M$, and set
$\psi(g)=\|i(gx)\|^2$. Then $\psi$ is conditionally negative
definite on $G$. By Proposition \ref{prop:resol_Trel}, $\psi$ is
bounded on $V$. This implies that the hypothesis of Theorem
\ref{thm:X^Q_non_vide} is fulfilled.\end{proof}

There are many metric spaces for which there automatically exists
such an equivariant embedding; namely, those metric spaces $M$
that have a $\text{Isom}(M)$-equivariant embedding in a Hilbert
space. Thus the hypotheses of Corollary are satisfied, for
instance when

$\bullet$ $M$ is a Hilbert space,

$\bullet$ $M$ is a tree, or a complete $\R$-tree \cite[Chap. 6,
Proposition 11]{HV}.

$\bullet$ $M$ is a real or complex hyperbolic space (maybe
infinite-dimensional) \cite{FH},

$\bullet$ $M$ is a finite-dimensional CAT(0) cube complex
\cite{NR}.

For instance, we have

\begin{cor}
Let $G,Q$ be locally compact groups, and let $f:G\to Q$ be a
resolution. Then $G$ has Property (FA) if and only if $Q$
does.\label{cor:resol_FA}
\end{cor}
\begin{proof} If $G$ has Property (FA), so does $Q$. Let us show the
converse. By a result of Alperin and Watatani (see \cite[Chap.
6]{HV}), every tree equivariantly embeds in a Hilbert space (more
precisely, the distance is a conditionally negative definite
kernel). It follows that, for every isometric action of $G$ on a
tree, there exists a nonempty $G$-invariant subtree on which the
action factors through $Q$. The result immediately
follows.\end{proof}

Theorem \ref{thm:resol} allows us to prove the converse of
Corollary \ref{cor:resol_irr}.

\begin{thm}
Let $f:G\to Q$ be a morphism between locally compact groups, with
dense image, and suppose $G$ $\sigma$-compact. Then $f$ is a
resolution if and only if, for every net $(\pi_i)$ of irreducible
unitary representations of $G$ converging to $1_G$, eventually
$\pi_i$ factors through a representation $\tilde{\pi}_i$ of $Q$,
and $\tilde{\pi}_i\to 1_Q$.\label{thm:resol_irr}
\end{thm}
\begin{proof} The condition is necessary by Corollary
\ref{cor:resol_irr}. Conversely suppose that it is
satisfied. Let us show that (2) of Theorem \ref{thm:resol} is
satisfied, using Theorem \ref{thm:T_rep_irr}. Fix $\eps>0$, let
$X\subset G$ be a subset such that $\overline{f(X)}$ is compact,
and let $(\pi_i)$ be a net of irreducible unitary representations
of $G$ converging to $1_G$. Then eventually $\pi_i$ factors
through a representation $\tilde{\pi}_i$ of $Q$, and
$\tilde{\pi}_i\to 1_Q$. Since $\overline{f(X)}$ is compact, this
implies that, eventually, $\tilde{\pi}_i$ has a
$(\overline{f(X)},\eps)$-invariant vector, so that $\pi_i$ has a
$(X,\eps)$-invariant vector.\end{proof}

%%%%%%%%%%%%%%%%%%%%%%%%%%%%%%%%%%%%%%%%%%%%%%%%%%%%%%%%%%%%%%%%%%%%%%%%%%%%%%%%%%

% ----------------------------------------------------------------
\end{document}